\documentclass[11pt]{article}
\usepackage{amsmath,amsfonts,amsthm,amssymb,mathtools, color}
\usepackage{enumerate}
\usepackage{tikz}

\usepackage{graphicx}

\def\ds{\displaystyle}
\def\forall{\hbox{for all}~}
\def\L{{\bf L}}

\def\ve{\varepsilon}

\def\bfP{{\bf P}}
\def\R{\mathbb{R}} 

\def\implies{\Longrightarrow}
\def\vp{\varphi}

\def\P{{\cal P}}
\def\caL{{\cal E}}
\def\TV{\hbox{Tot.Var.}}

\def\v{\vskip 1em}
\def\O{{\cal O}}
\def\D{{\cal D}}
\def\C{{\cal C}}

\def\T{{\cal T}}
\def\Q{{\mathcal Q}}
\def\ol{\overline}
\def\bega{\begin{array}}
\def\enda{\end{array}}
\def\begi{\begin{itemize}}
\def\endi{\end{itemize}}
\def\ov{\overline}
\def\dint{\int\!\!\!\!\int}

\def\Tilde{\widetilde}
\def\wt{\widetilde}
\def\Hat{\widehat}

\def\meas{\hbox{meas}}
\def\bel{\begin{equation}\label}
\def\eeq{\end{equation}}
\def\sqr#1#2{\vbox{\hrule height .#2pt
\hbox{\vrule width .#2pt height #1pt \kern #1pt
\vrule width .#2pt}\hrule height .#2pt }}
\def\square{\sqr74}
\def\endproof{\hphantom{MM}\hfill\llap{$\square$}\goodbreak}

\newtheorem{theorem}{Theorem}[section]

\newtheorem{lemma}[theorem]{Lemma}
\newtheorem{remark}[theorem]{Remark}
\newtheorem{prop}[theorem]{Proposition}
\newtheorem{example}[theorem]{Example}

\numberwithin{equation}{section}

\begin{document}
\title{\bf Intermediate Domains for Scalar Conservation Laws}

\author{Fabio Ancona$^{(1)}$, Alberto Bressan$^{(2)}$, Elio Marconi$^{(1)}$,  and
Luca Talamini$^{(1)}$\\~~\\
 {\small $^{(1)}$~Department of Mathematics 
 ``Tullio Levi-Civita", University of Padova,}\\
 {\small $^{(2)}$~Department of Mathematics, Penn State University,}\\
 \\~~\\
{\small E-mails: ancona@math.unipd.it,  ~axb62@psu.edu,}\\
{\small ~elio.marconi@unipd.it,~luca.talamini@math.unipd.it}\\~~~\\
}
\maketitle
\v

\begin{abstract} For a scalar conservation law with strictly convex flux, by  Oleinik's estimates
the total variation of a solution with initial data $\ov u\in \L^\infty(\R)$ 
decays like $t^{-1}$.  This paper introduces a class of intermediate domains 
$\cal P_\alpha$, $0<\alpha<1$, such that for $\ov u\in  \cal P_\alpha$ a faster decay rate is achieved: $\TV\bigl\{ u(t,\cdot)\bigr\}\sim t^{\alpha-1}$. A key ingredient of the 
analysis is a ``Fourier-type" decomposition of $\ov u$ into components which oscillate more and more rapidly.   The results aim at extending the theory of 
fractional domains for analytic semigroups to an entirely nonlinear setting.
\end{abstract}


\section{Introduction}
\label{sec:1}
\setcounter{equation}{0}
Consider a scalar conservation law 
\bel{claw}u_t+ f(u)_x~=~0,\eeq
 with strictly convex flux.
By a classical result \cite{Crandall, Kru}, there exists a contractive semigroup $S:\L^1(\R)\times \R_+
\mapsto \L^1(\R)$  such that, for every initial datum
\bel{idbu} u(0,\cdot) ~=~\ol u\in\L^1(\R),\eeq 
the
trajectory $t\mapsto u(t) = S_t \ol u$ is the unique entropy weak solution of the Cauchy problem.

It is well known that, even for smooth initial data, the solution can develop shocks in finite time.
Taking an abstract point of view, consider the operator $Au\doteq \frac{\partial}{\partial x} f(u)$
which generates the semigroup. 
Then there exists data
$\ol u\in Dom(A)$ in the domain of the generator, such that  $S_\tau\ol u\notin Dom(A)$
for some $\tau>0$.    In other words, the domain of the generator is not positively invariant.

To address this issue, 
the paper \cite{Cra2} introduced a definition of ``generalized domain" $\D$ for the operator.  This consists of all initial data $\ol u$ for which the map $t\mapsto S_t\ol u$
is globally Lipschitz continuous. Notice that for the conservation law (\ref{claw}) 
one has
$$\L^1\cap BV~\subseteq~\D.$$
Using the fact that the semigroup is contractive, it is easy to show that this generalized domain is 
positively invariant.   Indeed, the quantity
$$\limsup_{\ve\to 0+}~{\bigl\| S_{t+\ve}\ol u - S_t \ol u\bigr\|_{\L^1}\over \ve} $$
is a non-increasing function of $t$.

 Our present aim is to study intermediate domains
\bel{Da} \D_\alpha~\subset~\L^1(\R),\qquad\qquad  0<\alpha < 1,\eeq
related to the decay properties of the corresponding  trajectories of (\ref{claw}).
As in \cite{ABBCN} we define
\bel{Da2} \D_\alpha~\doteq~\Big\{ \ol u\in\L^1(\R)\,;~~\sup_{0<t<1} t^{-\alpha}  \bigl\| S_t \ol u-\ol u\|_{\L^1} < +\infty\Big\}.
\eeq
We also consider the domains (slightly different from the ones considered in \cite{ABBCN})
\bel{Da1} \Tilde\D_\alpha~\doteq~\Big\{ \ol u\in\L^1(\R)\cap\L^\infty(\R)\,;~~\sup_{0<t<1} t^{1-\alpha} \cdot \TV\bigl\{ S_t \ol u\} < +\infty\Big\}.
\eeq

The domains $\D_\alpha$  arise naturally in connection with balance laws:
\bel{nh1}
u_t+f(u)_x~=~g(t,x).\eeq
Indeed, as shown in \cite{ABBCN}, one has
\begin{prop} \label{l:1} Let $f\in \C^2$ with $f''(u)\geq c>0$ for all $u\in \R$.  
Consider a compactly supported solution $u=u(t,x)$ of (\ref{nh1}), and assume that
the source term
satisfies
\bel{gpb}\bigl\|g(t,\cdot)\bigr\|_{\L^1}~\leq ~C\qquad\qquad\forall t\in [0,T].\eeq   
Then for every $0<t\leq T$, one has $u(t,\cdot)\in \D_{1/2}$.
\end{prop}
We remark, however, that some of the functions $u(t,\cdot)$ can be unbounded.
In particular, they can have infinite total variation.

In the theory of linear analytic semigroups \cite{Henry, Lu}, intermediate domains
arise naturally as domains of fractional powers of sectorial operators. 
The faster decay of solutions is usually related to higher 
Sobolev regularity of the initial data.  In particular,  this theory applies to
semilinear equations of the form
\bel{CPlin}u_t - \Delta u~=~F(x,u,\nabla u),\qquad\qquad u(0,\cdot) = \ol u.\eeq
Under natural assumptions (see \cite{Henry} for details), this Cauchy problem is well posed provided that
the initial datum $\ol u$ lies in the domain of some fractional power 
 $(-\Delta)^\alpha$ of the generator.

Our eventual goal is to develop a similar theory of intermediate domains
for nonlinear semigroups generated by conservation laws.  
In particular, we conjecture that the Cauchy problem for a genuinely nonlinear
$2\times 2$ hyperbolic system with $\L^\infty$ initial data \cite{BCM, GL} 
is well posed within an intermediate domain such as (\ref{Da2})  or \eqref{Da1}. 

As a first step in this research program, here we focus our attention
on the scalar conservation law (\ref{claw}), seeking conditions on the initial data 
$\ol u\in \L^1(\R)$ that will imply $\ol u\in \D_\alpha$ or $\ol u\in \Tilde\D_\alpha$,
respectively.

Assumptions that imply $\ol u\in \D_\alpha$ can be readily
formulated in terms of fractional Sobolev regularity. 
On the other hand, conditions that guarantee a faster decay rate
of the total variation are more subtle.  Here we consider the assumption
\begi
\item[{\bf (P$_\alpha$)}]  {\it 
For every $\lambda\in\, ]0,1]$, there exists an open set $V(\lambda)\subset\R$ such that the following holds.
\bel{mop}
\meas\bigl(V(\lambda)\bigr) ~\leq~C\, \lambda^\alpha,
\eeq
\bel{tvp}
\TV\bigl\{ \ol u\,; ~\R\setminus V(\lambda)\bigr\}~\leq~C \,\lambda^{\alpha-1},\eeq
for some constant $C$ independent of $\lambda$.}
\endi
Roughly speaking, $\ol u$ can have unbounded variation, but most of its variation
should be concentrated on a set with small Lebesgue measure.
Our main result establishes the implication
\bel{impl}\ol u~~\hbox{satisfies}~{\bf (P_\alpha)}\qquad\implies\qquad \ol u\in \Tilde\D_\alpha\eeq
when $1/2<\alpha<1$.    On the other hand, a counterexample shows that the above implication
fails for $\alpha\leq 1/2$.
The proof of~\eqref{impl} 
relies on a structural result 
for functions satisfying 
${\bf(P_\alpha)}$, which is of independent interest. 
Indeed, Theorem~\ref{thm:dec}
provides a nonlinear ``Fourier-type" decomposition of such functions,
in components which oscillate more and more rapidly.

The remainder of the paper is organized as follows. In Section~\ref{s:mis}
we describe a general class of metric interpolation spaces, 
for functions defined on a set $\Omega\subseteq\R^N$.  This yields a 
natural way to formulate conditions such as {\bf (${\bf P}_\alpha$)}, in a general setting.

Section~\ref{s:3} contains some examples.  The first one (Fig.~\ref{f:ag34})
 shows how to construct an initial datum $\ol u$ with unbounded variation, such that $\ol u\in 
\Tilde \D_\alpha$, for any given $0<\alpha<1$.   We then consider initial data 
consisting of a packet of triangular waves (Fig.~\ref{f:ag47}).  By suitably 
choosing the size and distance of these triangular blocks we show that, 
if $0<\alpha\leq 1/2$, then there 
exists an initial datum
$\ol u$ that satisfies {\bf ($\bfP_\alpha$)} and yet $\ol u\notin \Tilde\D_\beta$
for any $\beta\in \,]0,1[\,$.
As stated in Proposition~\ref{p:33}, the implication (\ref{impl}) thus cannot hold for 
$\alpha\leq 1/2$.

Section~\ref{s:4}  is concerned with the intermediate domain $\D_\alpha$.   For $0<\alpha <1$ we
prove that any one of the conditions: (i) $\ol u$ lies in the fractional Sobolev space $W^{\alpha,1}(\R)$, or
(ii) $\ol u$  satisfies {\bf ($\bfP_\alpha$)},
implies $\ol u\in \D_\alpha$. These results are valid for any flux $f\in\C^1$, not necessarily convex.

Section~\ref{s:5} establishes further properties of functions
which satisfy {\bf ($\bfP_\alpha$)}, proving a  useful decomposition result, stated in
Theorem~\ref{thm:dec}.
Finally, in Section~\ref{s:6} we prove our main theorem, 
showing that for $1/2<\alpha<1$ 
the property {\bf ($\bfP_\alpha$)} implies that $\ol u\in \Tilde\D_\alpha$.
To simplify the exposition, the proofs will first be given for Burgers' equation.
In Remark~\ref{r:62} we observe that all results remain valid for a conservation law with
uniformly convex flux.

For an introduction to the theory of 
conservation laws we refer to \cite{Bbook, Btut, Daf, HR}.  Results on the decay of solutions to conservation laws in generalized BV spaces can be found in \cite{BGJ, Marconi}.
In addition to genuinely nonlinear conservation laws, several other examples
of nonlinear semigroups with regularizing properties are known in the literature,
see in particular~\cite{Barbu, BeCr, CP, Pa}.

\section{A family of metric interpolation spaces}
 \label{s:mis}
\setcounter{equation}{0}

Consider an open set $\Omega\subseteq\R^N$ and let $X$ be a Banach space contained in the set $L^0(\Omega)$ of Lebesgue measurable functions $f:\Omega\mapsto\R$. Let $0<\alpha <1$ be given.
A  distance function $d(\cdot, \cdot):L^0(\Omega)\times L^0(\Omega) \to [0,+\infty]$ 
can be defined as follows.  For any $\lambda\in \,]0,1]$,   we begin by 
setting
\bel{d5}
d^\lambda (f,g)~\doteq~d^\lambda (f-g,0),\eeq
\bel{d6}\bega{l} 
d^\lambda (f,0)~\doteq~\inf\Big\{ C\geq 0\,;~~\hbox{there exists $\Tilde f\in X$ such that}\\[3mm]
\qquad\qquad \qquad\qquad  \|\Tilde f\|_X\, \leq C\, \lambda^{\alpha-1}\,,\qquad \meas\bigl\{ x\in \Omega\,;~f(x)\not= \Tilde f(x)\bigr\}\, \leq\,
 C \, \lambda^{\alpha}\Big\}.\enda \eeq
Finally, we define
\bel{d7}
d(f,g)~\doteq~\sup_{0<\lambda \leq 1} ~d^\lambda(f,g).\eeq
By possibly identifying couples of functions $f,g,$ such that $d(f,g)=0$, we claim that (\ref{d7}) 
yields a distance on the set of Lebesgue measurable functions $f\in L^0(\Omega)$ for 
which $d(f,0)<\infty$.   This is usually
a strictly larger set than $X$.

\begin{lemma}   Let $0<\alpha<1$ be given, and let $d(\cdot, \cdot)$ be as in (\ref{d5})--(\ref{d7}).
Then the following properties hold.
\begi
\item[(i)]  $d(f,g)= d(g,f) \geq 0$.
\item[(ii)] If $f\in X$, then $d(f,0)\leq \|f\|_X$.
\item[(iii)] $d(f,h)~\leq~d(f,g) + d(g,h)$.
\endi
\end{lemma}

{\bf Proof.} {\bf 1.} Part (i) is trivial.  If $f\in X$, choosing $\Tilde f=f$ in (\ref{d6}), we see that $d^\lambda(f,0)\leq \|f\|_X$ for every $\lambda>0$. This yields (ii).
\v
{\bf 2.} We now check that each $d^\lambda(\cdot,\cdot)$, $0<\lambda\leq 1$, satisfies the triangle inequality.
Toward this goal,
let $f_1, f_2$ be measurable functions such that
$$d^\lambda(f_i,0) = C_i,\qquad i=1,2.$$
We then need to show that
\bel{tri}d^\lambda(f_1,f_2)~\doteq~d^\lambda(f_1-f_2,0)~\leq~C_1+C_2\,.\eeq
Given $\ve>0$, by assumption there exist 
functions $\Tilde f_1, \Tilde f_2\in X$ such that
\bel{f123}
\|\Tilde f_i\|_X\,\leq\,C_i \lambda^{\alpha-1}+\ve,
\qquad \meas\bigl\{ x\in \Omega\,;~\Tilde f_i(x)\not= f_i(x)\bigr\}\, \leq\, C_i \lambda^{\alpha}+\ve,\qquad \qquad i=1,2.\eeq
Then the function $\Tilde f_1 -\Tilde f_2$ satisfies
$$\|\Tilde f_1-\Tilde f_2\|_X~\leq~(C_1+C_2) \lambda^{\alpha-1}+2\ve,$$
$$ \meas\Big\{ x\in \Omega\,;~\Tilde f_1(x)-\Tilde f_2(x)\not= f_1(x)-f_2(x)
\Big\}~ \leq ~(C_1+C_2) \lambda^{\alpha} + 2\ve.$$
Since $\ve>0$ was arbitrary, this proves (\ref{tri}).
\v
{\bf 3.} 
In turn, the triangle inequality (iii)  follows from
$$\bega{rl} d(f-g,0)&=~\ds\sup_{\lambda \in ]0,1]} ~d^\lambda (f-g,0)~\leq~\sup_{\lambda \in ]0,1]} ~\Big( d^\lambda (f,0)+ d^\lambda (g,0)\Big)\\[4mm]
&\leq~\ds \sup_{\lambda \in ]0,1]} ~d^\lambda (f,0)+ \sup_{\lambda \in ]0,1]} ~d^\lambda (g,0)~=~d(f,0)+ d(g,0).\enda$$
\endproof

\begin{remark}\label{r:22}
 {\rm One should keep in mind that, in general, the  balls $\{ g\in {L^0(\Omega)}\,;~~d(g,f)\leq r\}$  are not convex. Moreover,  the function $f\mapsto d(f,0)$ is not a norm.}
\end{remark}

\v
In connection with  {\bf ($\bfP_\alpha$)}, for $0<\alpha<1$ 
we consider the distances $d^\lambda(f,g)$ as in (\ref{d5}), where now
\bel{d66}\bega{l} 
d^\lambda (f,0)~\doteq~\inf\Big\{ C\geq 0\,;~~\hbox{there exists $\Tilde f\in \L^1(\R)$ such that}\\[3mm]
\qquad\qquad \qquad\qquad  \TV \{\Tilde f\}\le C\, \lambda^{\alpha-1}, \quad \meas \{x \in \R\,;~f(x)\ne \Tilde f(x)\}\le C\, \lambda^\alpha\Big\}.\enda \eeq

Finally, given  $\ol u\in \mathbf L^1(\mathbb R)\cap \mathbf L^\infty(\mathbb R)$, we 
define
\bel{npa}
\Vert \ol u\Vert_{\P_{\alpha}}~\,\doteq~ \sup_{0<\lambda \leq 1} \,d^\lambda(\ol u, 0)
\eeq
and write $\ol u\in {\cal P}_{\alpha}$ if 
$\Vert \ol u\Vert_{{\cal P}_{\alpha}} < +\infty$.   Notice that this holds provided that 
$\ol u$ satisfies the
condition  {\bf ($\bfP_\alpha$)}.
Throughout the following, we shall use $\Vert \ol u\Vert_{{\cal P}_{\alpha}}$ as a convenient notation. However, as already pointed out
in Remark~\ref{r:22}, one should be aware that  $\|\cdot\|_{\P_\alpha}$  is not a norm.

\v

\section{Examples}\label{s:3}
\setcounter{equation}{0}
We present here some examples, to motivate the results proved in later sections.
We consider Burgers' equation
\bel{Bur}u_t + \left(u^2\over 2\right)_x~=~0.\eeq
Throughout the following, we use the semigroup notation $t\mapsto S_t \ov u$ 
to denote the solution of (\ref{Bur}) with initial data (\ref{idbu}).

\begin{figure}[ht]
  \centerline{\hbox{\includegraphics[width=8cm]{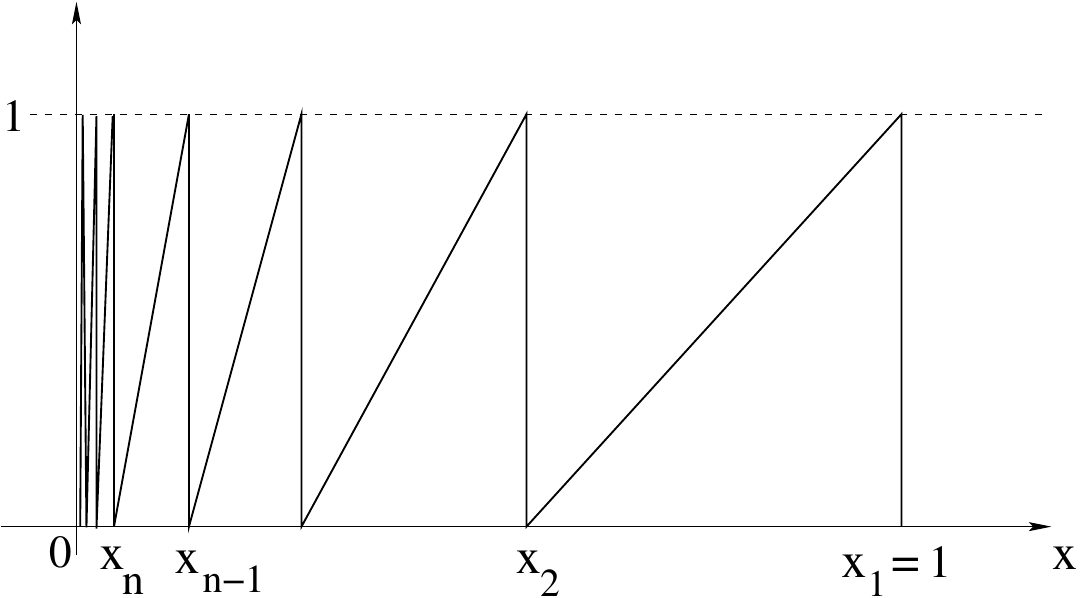}}}
  \caption{\small The initial data considered in Example~\ref{ex:1}.}\label{f:ag34}
\end{figure}

\begin{example}\label{ex:1} {\rm
Fix $\beta>0$ and consider the decreasing sequence of points $x_n= n^{-\beta} $, $n\geq 1$.
As shown in Fig.~\ref{f:ag34}, define the piecewise affine function
\bel{paf}
\ol u(x)~=~\left\{\bega{cl} 0\qquad &\hbox{if} \quad x\notin \,[0,1]\,,\\[4mm]
\ds {x-x_{n+1}\over x_n-x_{n+1}}\qquad &\hbox{if} \quad x_n<x<x_{n-1}\,.\enda\right.\eeq
We claim that this initial data lies in some of the subdomains $\Tilde \D_\alpha$,
depending on the exponent $\beta$.  Indeed, fix a time $t\in \,]0,1]$.
Consider the position $x_k(t)$ of the shock which is initially located at $x_k$.
By Oleinik's inequality, the total variation of the solution $u(t,\cdot)$ can be estimated by
\bel{oes1}\TV\bigl\{ u(t,\cdot)\,;~[0, x_k(t)]\bigr\} ~\leq~2{x_k(t)\over t}\,.\eeq
On the other hand, for $x>x_k(t)$, still by Oleinik's estimates we have
\bel{oes2}\TV\bigl\{ u(t,\cdot)\,;~[x_k(t), x_1(t) ]\bigr\} ~\leq~2{x_1(t)-x_k(t)\over (x_{k-1}-x_k)+ t}\,.\eeq
Observing that 
$$x_k(t)\,\leq\, x_k+t,\qquad\qquad x_1(t) - x_k(t)~\leq~1+t,$$
from (\ref{oes1})-(\ref{oes2}) we deduce
\bel{oes3} \TV\bigl\{ u(t,\cdot)\bigr\} ~\leq~2{x_k+t\over t} + 2{1+t\over (x_{k-1}-x_k)+ t}
~\leq~4 + 2{x_k\over t} + {2\over(x_{k-1}-x_k)+ t}\,. \eeq

Since we are assuming  $x_k=k^{-\beta}$,
the previous estimate yields
$$ \TV\bigl\{ u(t,\cdot)\bigr\} ~\leq~4 + {2k^{-\beta}\over t} + {2\over \beta k^{-\beta-1}+ t}\,.
$$
Here $k\geq 1$ is arbitrary.
Choosing $k\approx t^{-\gamma}$, we obtain
$$ \TV\bigl\{ u(t,\cdot)\bigr\} ~\leq~\O(1)\cdot \left( {t^{\beta\gamma}\over t} + {1\over t^{\gamma(\beta+1)}+ t}\right)~=~\O(1)\cdot \Big( t^{\beta\gamma-1} + 
t^{-\gamma(\beta+1)}\Big).$$
Here and throughout the sequel, the Landau symbol $\O(1)$ denotes a uniformly bounded quantity.
The two terms on the right hand side have similar magnitude if  $\gamma = (1+2\beta)^{-1}$.
With this choice, we obtain
$$ \TV\bigl\{ u(t,\cdot)\bigr\} ~\leq~\O(1)\cdot t^{-{\beta+1\over 2\beta+1}},$$
hence
$$\ol u\in \Tilde \D_\alpha,\qquad \hbox{with}\qquad \alpha~=~1 - {\beta+1\over 2\beta+1}~=~
{\beta\over 2\beta+1}.$$
}
\end{example}

\begin{figure}[ht]
  \centerline{\hbox{\includegraphics[width=14cm]{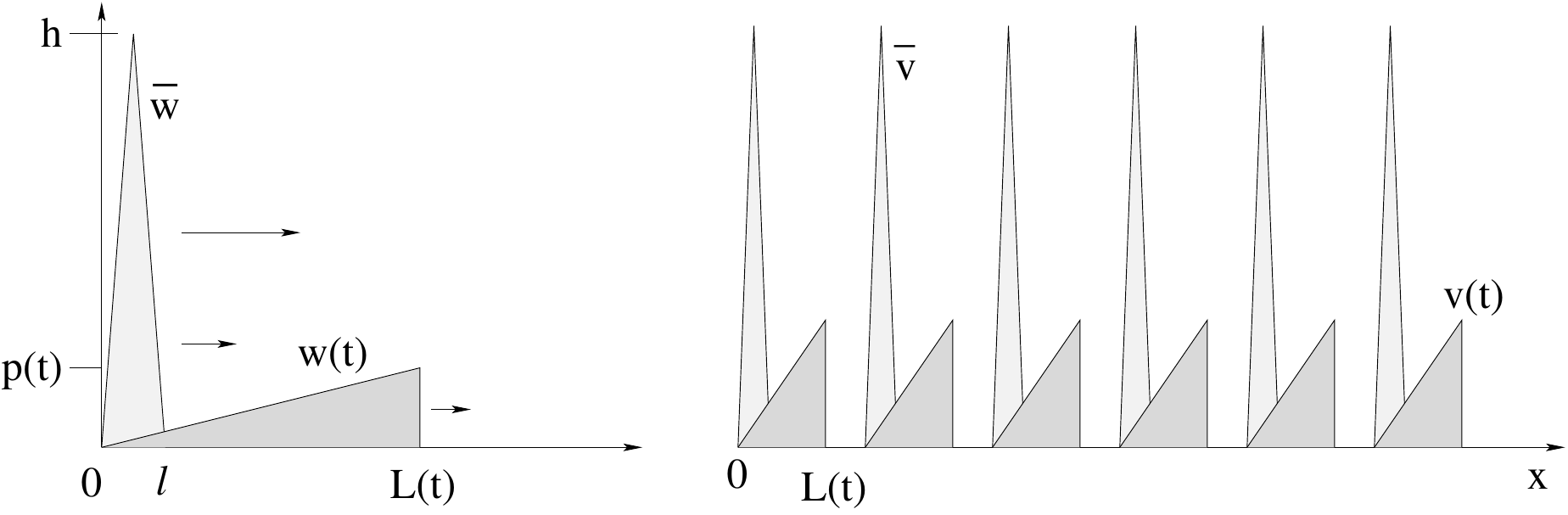}}}
  \caption{\small   Left: the elementary solution to Burgers' equation considered at (\ref{esol1})-(\ref{esol2}).  Right: the superposition of several shifted copies of the same solution.
}\label{f:ag46}
\end{figure}

In the next examples we consider initial data consisting of one or more triangular blocks.
As shown in Fig.~\ref{f:ag46}, left, the most elementary case is
\bel{esol1}
w(0,x)~=~\ov w(x)~=~\left\{ \bega{rl} h- \frac{2h}{\ell}\cdot |x-\ell/2| \quad &\hbox{if} ~~x\in [0,\ell],\\[1mm]
0\quad &\hbox{otherwise.}\enda\right.\eeq
At time $t=\ell/(2h)$, a shock is created in the solution at the point $\ell> 0$. Characteristics originating from points $0<x< \ell/2$  start impinging on the shock,
and the solution has a right triangle shape:
\bel{esol2} w(t,x)~=~\left\{ \bega{cl} \frac{2hx}{2ht+\ell}\quad &\hbox{if} ~~x\in [0,L(t)],\\[1mm]
0\quad& \hbox{otherwise}\enda\right.\eeq
Conservation of mass implies that the shock at time $t\geq \ell/(2h)$ is located at 
\bel{esol3} L(t)\,=\, \sqrt{ \frac{\ell}{2}\,(2\,h\,t + \ell)}.\eeq
Always for $t\ge\ell/(2h)$, we thus have 
 \bel{esol4} 
    \TV  \{S_t \ov w\} \,=\,  2\, p(t), \qquad\quad  p(t) \doteq h\sqrt{\frac{2\ell}{2ht+\ell}}.
 \eeq
 We notice for later use that the (decreasing) function $p(t)$ satisfies the lower bound
 \bel{esol5}
p(t) \, \geq\, \sqrt
    {h\, \ell\over  2t} \qquad \text{for all \; $t \geq \ell/2h$.}
 \eeq
 and that the (increasing) function $L(t)$ satisfies the upper bound
 \bel{esol6}
L(t) \, \leq \,\sqrt{2\,\ell\, h\, t} \qquad \text{for all \; $t \geq \ell/2h$.}
 \eeq
We also consider initial data containing packets of triangular blocks,
 shifted by different amounts so they do not overlap with each other. See Fig.~\ref{f:ag46}, right.

In the following proposition we denote by $\C^{0,\sigma}(\mathbb R)$ the space of Holder functions $\ol u: \mathbb R \to \mathbb R$ with exponent $0 < \sigma < 1$, equipped with the norm
$$
\Vert \ol u \Vert_{\C^{0,\sigma}} \doteq \Vert \ol u\Vert_{\C^0} + |\ol u|_{\C^{0,\sigma}}
$$
$$
|\ol u|_{\C^{0,\sigma}} \doteq \sup_{x < y}\frac{|\ol u(y)-\ol u(x)|}{|y-x|^{\sigma}}.
$$

\begin{prop}\label{p:33}
There exists a compactly supported function $\ol u: \mathbb R \to \mathbb R$ with the following properties:
\begin{enumerate}
    \item $\ol u \in  {\cal P}_{\alpha}$ for every $0 < \alpha \leq 1/2$;
    \item $\ol u \in \C^{0, \sigma}(\mathbb R)$ for every $0 < \sigma < 1$;
    \item $\ol u \notin \widetilde {\D}_{\beta}$ for any $0 <\beta < 1$. Namely:
    \bel{39}
    \limsup_{t \to 0+}  \;  t^{1-\beta} \cdot \TV \{S_t \ol u\}\, =\, +\infty \qquad \forall \; 0<\beta<1.
\eeq
\end{enumerate}
\end{prop}
\begin{remark}
    {\rm The function $\ol u$ constructed in Proposition \ref{p:33} does not belong to any $\P_{\alpha}$ if $1/2 < \alpha < 1$. This suggests that $\alpha = 1/2$ is a critical exponent for the decay of solution with initial data in $\P_{\alpha}$. In fact, this surprising behavior will be later confirmed by Theorem \ref{t:61}.}
\end{remark}
\begin{remark}
    {\rm By part 2.~of Proposition~\ref{p:33} and by the embedding $\C^{0,\sigma}(\mathbb R) \hookrightarrow W^{s, p}_{loc}(\mathbb R)$ for every $0 < s < \sigma < 1$, $p \geq 1$, we obtain that 
    $$
W^{s,p}(\mathbb R) \not\subset \widetilde {\D}_{\beta} \qquad \text{for every $0< s <1$, ~$1 \leq p < \infty$ and $0 < \beta <1$.} 
    $$
   On the other hand, the inclusion $W^{\alpha, 1}(\mathbb R) \subset \D_{\alpha}$ does hold for every $0 < \alpha < 1$, as proved in Proposition \ref{p:42}.}
\end{remark}

The proof of Proposition \ref{p:33} is based on the following lemma.
\begin{lemma}\label{l:34}
    For every fixed $t \in (0, 1)$, there exists a function $\widehat u: \mathbb R \to \mathbb R$ satisfying
    \bel{pab8} \left\{\bega{cl} \widehat u(x) \in [0,1]\quad &\hbox{if} ~~x\in [0,1],\\[1mm]
    \widehat u(x) =0 \quad &\hbox{if} ~~x\notin [0,1],\enda\right.\eeq
  \bel{pab4}\Vert \widehat u\Vert_{{\cal P}_{\alpha}}~ \leq~ C_0\,,\qquad \forall \; 0 < \alpha \leq 1/2, \qquad \qquad
    \TV \{S_t \widehat u \}~\geq~ \frac{1}{C_0\, t}\,.
  \eeq
  \bel{pab9}
 \|\widehat u\|_{\C^{0, \sigma}}~ \leq ~C_\sigma\, \doteq\, \left(e^{ 2^{\left(\frac{1}{1-\sigma}\right)}}\right)^{(1-\sigma) e^{-1}} \qquad \text{for all }\sigma \in \,]0,1[.
  \eeq
  Here $C_0$ is a constant independent of $t$.
\end{lemma}

We postpone the proof of the lemma, and begin by showing that it implies
the previous proposition.

{\bf Proof of Proposition~\ref{p:33}.} 
Let $(t_j)_{j\geq 1}$ be a sequence decreasing to $0$ sufficiently fast (to be specified later). Let
$$
x_0 = 0, \quad x_j \doteq \sum_{k=0}^{j-1} 2 \cdot 2^{-k}, \qquad I_j \doteq [x_j, x_{j+1}], \qquad j = 0,1,2 \ldots
$$
Let $\widehat u_{j}$ be a function satisfying the properties (\ref{pab8})-(\ref{pab4}) in Lemma~\ref{l:34} for $t = t_j$.   Consider the rescaled functions
$$
\ol u_{j}(x)~ \doteq ~2^{-j} \widehat u_{j}\big(2^{j}(x-x_j)\big), \qquad  j = 0,1,2 \ldots
$$
By Lemma~\ref{l:34}, the corresponding rescaled solution of Burgers' equation
satisfies 
$$\TV\{  S_{t_j} \ol u_{j} \}~=~2^{-j}\,\TV\{ S_{t_j} \Hat u_{j}\}~\geq~ \frac{1}{C_0\,2^j \, t_j}\,.$$ 
Define the initial data
\begin{equation}\label{e:olu}
\ol u ~= ~\sum_{j=0}^{\infty} \ol u_{j}.
\end{equation}
Notice that 
$${\rm supp} \; \ol u_{j}(t,\cdot)~ \subseteq~ [x_j, x_j + 2^{-j}(1+t)  ] ~\subset~ I_j
\qquad \qquad\hbox{for }~t\in [0,1].$$
In particular for every $i\ne j$ and $t< 1$ the supports of $S_t \ol u_{i}$ and $S_t\ol u_{j}$ remain disjoint.
Choosing $t_j < 2^{-j}$ we thus obtain
\[
 \bigl(S_{t_j} \ol u\bigr) (x) ~=~\bigl(  S_{t_j} \ol u_{j} \bigr)(x) \qquad \forall x \in I_j\,.
\]
Since every $\Hat u_j$ satisfies (\ref{pab4}), for $0 < \alpha \leq 1/2$, by the triangle inequality it now follows
$$
\Vert \ol u\Vert_{{\cal P}_{\alpha}} ~\leq~ \sum_{j=0}^{\infty} \Vert \ol u_{j}\Vert_{{\cal P}_{\alpha}}~ \leq ~C_0\sum_{j=0}^{\infty} 2^{-j} ~<~ +\infty.
$$
This implies that $\ol u$ satisfies assumption 1.~of Proposition~\ref{p:33}. By \eqref{pab9} of Lemma~\ref{l:34}, we infer that 
$$\|\ov u_j\|_{\C^{0,\sigma}}~ \leq ~2^{-(1-\sigma) j}\, \|\widehat u_j\|_{\C^{0,\sigma}} ~\leq~ C_\sigma \,2^{-(1-\sigma) j}.$$
Therefore since the supports of $\ol u_j$ are all disjoint, the function $\ol u$ belongs to all the H\"older spaces $\C^{0,\sigma}$ for $0 < \sigma < 1$.
Finally, choosing for example $t_j = \exp (-2^{j})$, for any $0<\beta<1$ we obtain
$$\bega{l}\ds
\limsup_{t \to 0+} \,t^{\beta}\,\TV\{   S_t\ol u\} ~ \geq ~\lim_{j \to +\infty}\, t_j^{\beta} \,\TV\{ S_{t_j} \ol u\} ~\geq  ~ \lim_{j \to +\infty} \,t_j^{\beta}\, \TV\{ S_{t_j} \ol u_j\} \\[4mm]
\ds\qquad  \geq~  \frac{1}{C_0\, 2^j \,t_j^{1-\beta}}~ \to~ +\infty,\enda
$$
completing the proof of 3. of Proposition~\ref{p:33}. \endproof

{\bf Proof of Lemma~\ref{l:34}.}

\begin{figure}[ht]
  \centerline{\hbox{\includegraphics[width=14cm]{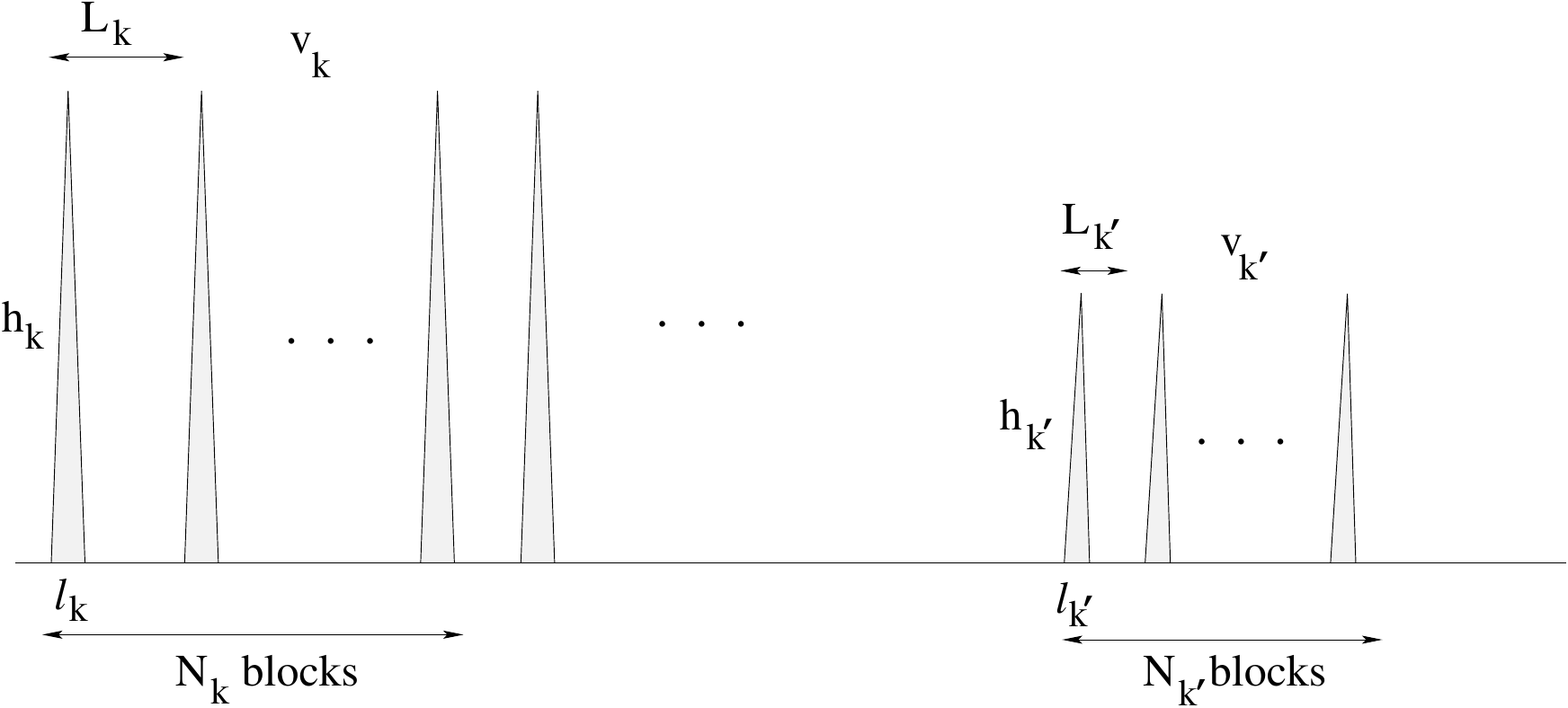}}}
  \caption{\small  A family of triangular wave packets.
}\label{f:ag47}
\end{figure}

{\bf 1.} Let $t > 0$ be a fixed positive time. Given a sequence of positive numbers $\{\ell_k\}_k$ satisfying $\ell_k \leq 2^{-k}$ (to be chosen later), and an integer $k \geq 1$, we  construct a packet of triangular waves by setting 
   \bel{vkw}
    v_k~ = ~\sum_{j=1}^{N_k} w_k^j\,,
    \eeq
    where $w_k^j = \ol w_k(x-j\,L_k)$ are translations by $j\,L_k> 0$ of elementary triangular blocks as in (\ref{esol1}), with width $\ell_k$ and height 
     \begin{equation}\label{eq:boxk}
    h_k \doteq 2^k \ell_k \leq 1\,.
    \end{equation}
    The distance $L_k$ between the supports of two blocks is  chosen large enough so that the supports of the corresponding solutions remain disjoint up to the given time $t$. By \eqref{esol6},  it is sufficient to separate these elementary blocks by a distance $L_k \doteq \sqrt{2 h_k \, \ell_k \, t}$, see Figure~\ref{f:ag47}.  With this choice, the support of $v_k$ 
    is contained inside an interval 
   $I_k$ with length      
   \bel{pak}
    \meas(I_k)  ~\leq ~ N_k L_k~ = ~\sqrt{2t} \,2^{k/2}N_k\ell_k\,.
    \eeq
Moreover, since the elementary blocks do not interact with each other up to time $t$, 
assuming $\ell_k/(2h_k) \leq t$, by \eqref{esol5}, one has
\bel{eq:TVvk}    \TV\{  S_t v_k \} ~=~ 2 N_k p_k(t) ~\ge~ \sqrt{2} N_k\cdot \sqrt{\frac{h_k \ell_k}{t}}
~=~
\frac{\sqrt{2}\cdot 2^{k/2}\,\ell_k}{\sqrt{t}}\, N_k\,.
\eeq
We now choose
\bel{Nk}
\ell_k \doteq 2^{-1} \cdot 2^{-k/2} \cdot k^{-k}, \qquad \qquad 
N_k \doteq 2^{-1} \cdot 2^{-k/2}\cdot \ell_k^{-1} = k^k
\eeq
and  notice that, by \eqref{eq:boxk}, this implies 
$$
 \TV\{ v_k\}~=~ 2 N_k h_k  ~=~  2^{k/2}.
$$

 \v
 {\bf 2.}
We put next to each other all the wave packets $v_k$, for $k_1 \leq k \leq k_2$, where $k_1, k_2 \in \mathbb N$ will be chosen later (see Figure \ref{f:ag47}), and define
\bel{buvk}
\widehat u ~\doteq ~\sum_{k=k_1}^{k_2} v_k\,.
\eeq
Recalling (\ref{d66})-(\ref{npa}), we now show that, for $\alpha=\frac12$, the above construction yields a uniform bound
\bel{upab}\|\widehat u\|_{\P_\alpha}~\leq~C,\eeq
with a constant $C$ independent of $k_1,k_2$. As a consequence, the same bound holds for $\alpha \in \, ]0,\frac12[$.

To prove (\ref{upab}),
let $\lambda \in \, ]0,1]$ and choose $\bar k$ such that $\lambda \in \, ]2^{-\bar k}, 2^{-\bar k +1}]$.
Set $V(\lambda)$ in \eqref{mop}-\eqref{tvp} to be the support of $\sum_{k=\bar k}^{k_2} v_k$. By \eqref{Nk} the following estimates hold:
\[
\meas(V(\lambda))~ \le ~\sum_{k=\bar k}^{k_2}N_k \ell_k ~=~ \frac{1}{2}\sum_{k=\bar k}^{k_2}2^{-k/2} ~<~ \frac{1}{2}\frac{\sqrt 2}{\sqrt 2 -1}\cdot 2^{-\bar k/2} ~<~ C_0\, \lambda^{1/2},
\]
\[
\TV\left\{ \sum_{k=k_1}^{\bar k-1} v_k;\, \R\right\}~ \le~ \sum_{k=k_1}^{\bar k-1}N_k 2^k \ell_k ~<~ \frac{1}{2}\frac{2^{\bar k/2}}{\sqrt 2 -1}~<~ C_0 \,\lambda^{1/2}.
\]
This proves (\ref{upab}).
A more general result will be obtained in Proposition \ref{prop:dec}, to which we refer for additional details.
\v
{\bf 3.}
In view of (\ref{pak}), \eqref{Nk}, the support of $\widehat u$ is contained in an interval $I$ whose
length is 
\begin{equation}\label{eq:space}
\meas(I)~\leq~\sum_{k_1}^{k_2} \meas(I_k)~=~
\sum_{k_1}^{k_2} \sqrt{2t} 2^{k/2}N_k\ell_k ~=~ \sqrt{\frac{t}{2}}(k_2-k_1+1).
\end{equation}
We now compute the total variation of the corresponding solution 
$S_t \Hat u$ at a given time $t>0$.   ~ If $t\ge 2^{-1}\cdot 2^{-k_1}$, i.e.,  if $k_1 \geq \log_2(1/2t)$, then at time $t$ all the elementary solutions appearing in the blocks $v_k$, $k \geq k_1$, have a right triangle shape and we can use the estimate \eqref{eq:TVvk}.  Since these blocks do not interact with each other, we have
\bel{eq:TVt}
\TV\{  S_t \widehat u\}~ =~ \sum_{k_1}^{k_2} \TV  \{S_t v_k \}~\ge~ \sqrt{\frac2t}\sum_{k_1}^{k_2} 2^{k/2}N_k \ell_k ~=~ \sqrt{\frac{1}{2t}}(k_2-k_1+1).
\end{equation}
Using the notation  $\lceil a \rceil$ to denote the smallest integer $\geq a$,
we now choose $k_1= \lceil \log_2(1/t)\rceil$ and $k_2= k_1 -2 +\lceil\sqrt{\frac{2}{t}}\rceil$. 
By \eqref{eq:space} we deduce that the support of $\widehat u$ is contained in the interval $I$ whose length is
$$
\mathrm{meas}(I) \leq \sqrt{\frac{t}{2}}(k_2-k_1+1) = \sqrt{\frac{t}{2}}\left( \left\lceil\sqrt{\frac{2}{t}}\right\rceil-1\right) \leq 1.
$$
By \eqref{eq:TVt}, $S_t \widehat u$ has total variation 
$$
\TV\{  S_t \widehat u\}~ \geq  \sqrt{\frac{1}{2t}}(k_2-k_1+1) = \sqrt{\frac{1}{2t}}\left( \left\lceil\sqrt{\frac{2}{t}}\right\rceil-1\right) \geq \sqrt{\frac{1}{2t}}\left\lceil \frac{\sqrt{2}-1}{\sqrt{t}} \right\rceil \geq \frac{1}{c} \cdot\frac{1}{t}
$$
where $c> 0$ is an absolute constant. 

Finally, by \eqref{eq:boxk}, \eqref{Nk}, 
the function $\widehat u$ satisfies 
$$
\sup_{x < y}\frac{|\widehat u(y)-\widehat u(x)|}{|y-x|^{\sigma}} ~\leq~ \sup_{k} \big( k^{-(1-\sigma)k}\cdot 2^{k}\big) ~\leq~\big(e^{ 2^{\frac{1}{1-\sigma}}}\big)^{(1-\sigma) e^{-1}},
$$
since the function $k \mapsto k^{-(1-\sigma)k}\cdot 2^k$ attains its maximum in $\R$ at $k=2^{\frac{1}{1-\sigma}} e^{-1}$. 
\endproof

\section{The intermediate domains $\D_\alpha$}
\label{s:4}
{
We consider entropy solutions to the scalar conservation law \eqref{claw}}
and let $S: \L^{1}(\mathbb R) \times \R_+\to \L^1(\R)$
be the corresponding semigroup of entropy weak 
solutions.
{
Notice that in this section the convexity assumption of the flux $f$ is not necessary.}

The goal of this section is to study the subdomains ${\cal D}_{\alpha} \subset \mathbf L^1(\mathbb R)\cap \mathbf L^{\infty}(\mathbb R)$ defined by 
\bel{Dali}
    {\cal D}_{\alpha} \doteq \Bigg\{ \ol u \in\mathbf L^1(\mathbb R)\cap \mathbf L^{\infty}(\mathbb R) \; ;~~   \sup_{0 < t \leq 1} t^{-\alpha} \Vert  S_t \ol u-\ol u\Vert_{\mathbf L^1} < +  \infty \Bigg\}.
\eeq

Since the semigroup $ S_t$ is nonlinear, the domains ${\cal D}_{\alpha}$ are not 
vector spaces.
However, we can ask if they contain some classical linear spaces, such as fractional Sobolev spaces. 
\smallskip

Let  $1\leq p < +\infty$ and $0<\alpha \leq 1$ be given, together with an open set 
$\Omega\subset\R$.  The fractional Sobolev space $W^{\alpha, p}(\Omega)$ 
is defined by  (see for example \cite{DPV_sobolev})
\begin{equation}
    W^{\alpha, p}(\Omega) ~\doteq~
     \Bigg\{u \in \L^p(\Omega) ~ ;~ ~ \frac{|u(x)-u(y)|}{|x-y|^{\frac{1}{p}+\alpha}}\in \L^p(\Omega\times \Omega)   \Bigg\},
\end{equation}
equipped with the norm
\begin{equation}
    \Vert u\Vert_{W^{\alpha, p}} ~\doteq~ \|u\|_{\L^p}+   \Bigg(\int_{\Omega \times \Omega} \frac{|u(x)-u(y)|^p}{|x-y|^{1+\alpha p}}\, dxdy \Bigg)^{\frac{1}{p}}.
\end{equation}

As it is well known, functions in Sobolev spaces can be approximated by smooth functions by taking mollifications.
Let $\eta:\R\mapsto [0,1]$ be a symmetric, $\C^\infty$ mollifier 
with compact support, so that
\bel{moll}\left\{ \bega{cl}\eta(s)= \eta(-s)\,\in\,[0,1]\qquad &\hbox{if} \quad 
s\in [-1,1], \\[2mm]
\eta(s)\,=\,0\qquad &\hbox{if}\quad |s|\geq 1,\\[2mm]
\bigl|\eta'(s)\bigr|\leq 2\qquad &\forall s\in\R,
\enda\right.\qquad\qquad \int\eta(s)\, ds=1.\eeq
Here and in the sequel, the prime $'$ denotes a derivative.
For $h>0$, define the rescaled kernels by setting
\bel{resk} \eta_h(s)~=~{1\over h} \eta\left({s\over h}\right).\eeq
For $u\in \L^1_{loc}$, consider the convolution $u_h = u\star \eta_h$.
The rate of
convergence of these mollifications depends on the regularity properties of the 
function $u$. 
\begin{lemma}\label{l:51}
Assume $ u \in  W^{\alpha, 1}(\R)$ for some $0<\alpha \leq 1$. Then, for every $h > 0$,  the convolution $u_h = u\star \eta_h$ satisfies
\begin{equation}\label{eq:approquant}
    \Vert  u- u_{h}\Vert_{\L^1} ~\leq~ \Vert u\Vert_{W^{\alpha, 1}} \cdot h^{\alpha}, \qquad \qquad   \Vert  u_{h}^{\prime}\Vert_{\L^1} ~ \leq~ C \Vert u\Vert_{W^{\alpha, 1}}\cdot  \frac{1}{h^{1-\alpha}},
\end{equation}
for some constant $C$ independent of $u$.
\end{lemma}

{\bf Proof.}
A direct computation yields 
\[
        \bega{rl}\ds
            \int |u(x) - u_{h}(x) |\,dx &\ds \leq ~\dint |u(x)-u(y)| \eta_{h}(y-x) \,dxdy \\[4mm]
            & = ~\ds\int \frac{1}{h}\eta(s/h)  \int |u(x+s)-u(x)| \,dx ds \\[4mm]
            &\ds \leq ~\int_{{-h}}^h \left( \frac{1}{|s|} \int \bigl|u(x+s)-u(x)\bigr| \,dx
            \right)ds \\[4mm]
            &\ds \leq ~h^{\alpha} \dint \frac{|u(x+s)-u(x)|}{|s|^{1+\alpha}}  \,dxds \\[4mm]
            &=~ \Vert u\Vert_{W^{\alpha, 1}} \cdot h^{\alpha}.
        \enda \]
Moreover, the total variation of $u_h$ is bounded by
\[
    \bega{rl}\ds
        \int_\R |u_h^{\prime}(x) | \,dx &\ds =~ \int \left|\int  u(y) \eta^{\prime}_h(y-x) dy  
        \right| \,dx 
        \\[4mm]
        & =~\ds
        \int \Big|\int u(x+s) \eta^{\prime}_h(s) ds   \Big| dx
        \\[4mm]
        &\ds =~ \int
         \frac{1}{h^2} \Big|   \int u(x+s)  \eta^{\prime}(s/h) ds \Big| dx\\[4mm]
        & = ~\ds 
        \int \frac{1}{h^2} \Big| \int_0^h \eta^{\prime}(s/h) (u(x+s)-u(x-s)) ds \Big| dx\\[4mm]
        & \leq ~\ds C \, \int \frac{1}{h^{1-\alpha}} \int_0^h  \frac{|u(x+s)-u(x-s)|}{s^{1+\alpha}} ds dx \\[4mm]
        & \ds\leq~C \,\Vert u\Vert_{W^{\alpha, 1}} \frac{1}{h^{1-\alpha}},
    \enda
\]
Notice that the constant $C$ depends only on the mollifying kernel $\eta$.
\endproof

\v

\begin{prop}\label{p:42} Let (\ref{claw}) be any conservation law with continuously differentiable flux.
For every $\alpha \in \,]0, 1]$ we have the inclusion 
$ \L^{\infty}(\R) \cap W^{\alpha, 1}(\R) \subseteq {\cal D}_{\alpha}$. 
\end{prop}
{\bf Proof.}  Let $\ol u\in \L^{\infty}(\R) \cap W^{\alpha, 1}(\R) $ and consider the 
mollifications 
$u_h\doteq \eta_h\star \ol u$. By Lemma~\ref{l:51} it follows
\[
\Vert   S_h  u_h-u_h\Vert_{\L^1}~ \le ~\| f'(u_h)\|_{\mathbf L^\infty}\, h \cdot \TV \{u_h\}~ \le~ C\,\| f'\|_{\L^\infty} \| \ol u\|_{W^{\alpha,1}} h^\alpha,
\]
where the $\mathbf L^\infty$ norm of $f'$ is taken on the interval $[-\| \ol u\|_{\L^\infty},\| \ol u\|_{\L^\infty}].$
Therefore
\[\bega{rl}
 \Vert S_h \ol u-\ol u\Vert_{\mathbf L^1}& \leq ~ \Vert   S_h \ol u-   S_h  u_h\Vert_{\L^1} +\Vert   S_h  u_h-u_h\Vert_{\mathbf L^1}+\Vert u_h - \ol u\Vert_{\L^1} \\[2mm]
 &\leq ~\left( 2 +C\,\| f'\|_{\L^\infty} \right) \,\| \ol u\|_{W^{\alpha,1}} h^{\alpha}. 
 \enda
\]
\endproof
%
%
%
%
%
%
%
%
%
\v
The second result in this section is formulated in terms of the property 
{\bf (P$_\alpha$)}.

\begin{prop}\label{p:43} Let (\ref{claw}) be a conservation law with continuously differentiable flux.
For any $0<\alpha<1$, if $\ol u\in \L^\infty(\R)$ satisfies  {\bf (P$_\alpha$)}, then $\ol u\in \D_\alpha$.
\end{prop}

{\bf Proof.}  Let $\ol u$ satisfy {\bf (P$_\alpha$)}.   Given $t\in \,]0,1]$, set
$\lambda=t$ and let $V(\lambda)\subset\R$  be an open set satisfying (\ref{mop})-(\ref{tvp}).
%
Observing that this open set $V(t)$ is a countable union of disjoint open intervals
$$V(t)~=~\bigcup_{k\geq 1} \, ]a_j, b_k[ \,,$$
we  define a new function $\ol v$ by replacing $\ol u$ with an affine function on each  interval $[a_j, b_j]$.   Namely, 
$$\ol v(x)~=~\left\{ \bega{cl} \ol u(x)\quad &\hbox{if}\quad x\notin \cup_k [a_k, b_k],\\[1mm]
\ds {(b_j-x) \ol u(a_j) + (x-a_j) \ol u(b_j)\over b_j-a_j}\quad &\hbox{if}\quad  x\in [a_j, b_j]\,.\enda\right.$$
This implies
$$\TV\{\ol v\}~\leq~\TV\bigl\{ \ol u\,; ~\R\setminus V(t)\bigr\}~\leq~C\, t^{\alpha-1},$$
$$\|\ol v-\ol u\|_{\L^1}~\leq~ 2 \|\ol u\|_{\L^\infty} \cdot  \meas(V(t))~\leq~2 \|\ol u\|_{\L^\infty} \cdot C t^\alpha.$$

We thus obtain
$$\bega{rl} \| S_t \ol u - \ol u\|_{\L^1}&\leq~\| S_t \ol u - S_t\ol v\|_{\L^1}+ \|S_t \ol v - \ol v\|_{\L^1}
+\|\ol u - \ol v\|_{\L^1}\\[2mm]~&\leq~ \|\ol v-\ov u\|_{\L^1} + t\cdot \|f'\|_{\L^\infty} \cdot\TV\{\ol v\} +\|\ol v-\ov u\|_{\L^1} \\[2mm]
& \leq~4 \, \|\ov u\|_{\L^\infty} \, C\, t^\alpha + \|f'\|_{\L^\infty}\, C t^\alpha,
\enda$$
where the $\mathbf L^\infty$ norm of $f'$ is taken on the interval $[-\| \ol u\|_{\L^\infty},\| \ol u\|_{\L^\infty}]$.
Since the same constant $C$ is valid for all $t\in \,]0, 1]$, this proves 
that $\ol u\in \D_\alpha$.
\endproof
\v

\section{A decomposition property for functions  $\ol u\in {\cal P}_{\alpha}$}
\label{s:5}
In this section we study properties of functions that lie in the metric space $\P_\alpha$
introduced at (\ref{npa}). These are functions that satisfy the property {\bf ($\bfP_\alpha$)} at (\ref{mop})-(\ref{tvp}).  Our main result provides a decomposition
 of a function $\ol u \in \P_{\alpha}$, as the sum of countably many components with 
 different degrees of regularity.
 
\begin{theorem}\label{thm:dec}
     Let $\bar u:\R\mapsto\R$ be a  measurable function and let $0<\alpha<1$ be given.
Then  $\ol u\in {\cal P}_{\alpha}$ if and only if it can be decomposed as 
    \begin{equation}\label{udec}
    \ol u(x) ~= ~\sum_{k=0}^{\infty} v_k(x) \qquad \text{for a.e. $x \in \mathbb R$,}
    \end{equation}
    where the $v_k$ satisfy the following properties.  For some constant $C=\O(1)\cdot \Vert \ol u \Vert_{{\cal P}_\alpha}$ one has
    \begi
        \item[{\bf (i)}] {\bf - bounds on the support and on the total variation:} 
         \bel{v00}
    \mathrm{Tot.Var.}\{v_0\}~ \leq ~ C,
    \eeq
     \begin{equation}\label{vk}
      \mathrm{Tot.Var.}\{v_k\}\, \le\,  C\cdot 2^{(1-\alpha) k}, \qquad \hbox{\rm meas}\bigl(\{v_k \ne 0\}\bigr) \,\leq\, C \cdot  2^{-\alpha k}, \qquad \forall \; k \geq 1.
    \end{equation}
    \item[{\bf (ii)}] {\bf - one-sided Lipschitz bound:} 
    \bel{oslip}
v_k(x_2)-v_k(x_1)~\leq~2^k \cdot  (x_2-x_1)\qquad\quad\forall x_1<x_2\,.\eeq
  
   \item[{\bf (iii)}] {\bf -  a further decomposition:} 

         For each $k\ge 1$ we
 can further decompose
    $$
     v_k \,=\, \sum_{p=1}^{\infty} v_k^p
    $$
    so that the following conditions hold:  the functions $v_k^p$ satisfy \eqref{oslip}, their supports have disjoint interiors, and setting
        $
        \ell_k^p \doteq \meas({\rm supp}\, v_k^p)
        $
         it holds
        $$
        \bigl|v_k^p(x)\bigr|~ \leq~ h_k^p~ \doteq ~  2^k \ell_k^p, \qquad \forall \; x \in {\rm supp}\, v_k^p\,,
        $$
        and
        \begin{equation}\label{eq:rectsum}
        \sum_{p \geq 1} \ell_k^p~ \leq~ C \cdot 2^{-\alpha k}.
        \end{equation}
          \endi
\end{theorem}

\medskip

The proof of Theorem \ref{thm:dec} will be achieved in three steps. We first show in Lemma \ref{prop:dec} that the existence of a decomposition as in \eqref{udec} which satisfies property {\bf (i)}  is equivalent to the statement that 
$\ol u\in \P_{\alpha}$. Next, in Lemma \ref{lemma:strong_dec}  we show that this decomposition can be refined so to satisfy also property {\bf (ii)}, 
still with some constant $C$ of the same order of $\Vert \ol u\Vert_{\P_{\alpha}}$. Finally, Lemma \ref{l:54} shows that {\bf (iii)}  is an easy consequence of  {\bf (i)}  and {\bf (ii)}  

\medskip

   \begin{remark}\label{rem:embedding}{\rm
  One can estimate the $\mathbf L^p$ norm of $\ol u \in {\cal P}_{\alpha}$,  for 
  $1\leq p<+\infty$, by 
$$
    \Vert \ol u\Vert_{\L^p}~ \leq ~\sum_{k=0}^{+\infty}\Vert v_k\Vert_{\L^p}~ 
    \leq ~\sum_{k=0}^{+\infty}\Vert v_k\Vert_{\L^\infty} \meas(\{ v_k \ne 0\})^{\frac1p} \leq ~ \frac{C}2 \cdot C^{1/p}\cdot \sum_{k=0}^{+\infty}  2^{k(1-\alpha -\alpha p^{-1})},
    $$
    where $C=\O(1)\cdot\Vert \ol u\Vert_{{\cal P}_{\alpha}}$ is the same constant as in Theorem \ref{thm:dec}.
 If $p < \frac{\alpha}{1-\alpha}$, we thus have the embedding ${\cal P}_{\alpha} \hookrightarrow \L^p_{loc}$.   Indeed, for every compact set $K \subset \R$, 
 there holds
$$
\Vert \ol u\Vert_{\L^p(K)}~ \leq~ c(K) \cdot \frac{1}{1-2^{(1-\alpha-\alpha p^{-1})}}\cdot \Vert \ol u\Vert_{{\cal P}_{\alpha}}^{1 + \frac1p}\,.$$

In particular ${\cal P}_{\alpha} \hookrightarrow \mathbf L^1_{loc}$ if  $\alpha > 1/2$. 
Notice that this is consistent with the scaling property
\[
\Vert u \Vert_{{\cal P}_{\alpha}}\,=\, \Vert u_\mu \Vert_{{\cal P}_{\alpha}} \qquad \mbox{with} \quad u_\mu(x)\doteq \mu^{\frac{1-\alpha}\alpha}u(\mu x)\qquad  \text{for $\mu \geq 1$}.
\]
More generally, if $p < \frac{\alpha}{1-\alpha}$, the immersion ${\cal P}_{\alpha} \hookrightarrow \mathbf L^p_{loc}$ is compact, namely a bounded sequence $\{\ol u_n\}_{n \in \mathbb N} \subset {\cal P}_\alpha$ admits a convergent subsequence in $\mathbf L^p_{loc}$.

Indeed, if $\ol u_n$ is a sequence of functions with $\Vert \ol u_n\Vert_{{\cal P}_{\alpha}}$  uniformly bounded in $n$, and if $\{v_k^n\}_{k \in \mathbb Z}$ are the functions appearing in the decomposition of $\ol u_n$, by a diagonal argument using Helly's compactness theorem one extracts a subsequence $\{n_i\}_{i \in \mathbb N}$ such that $v_k^{n_i}$ converges in $\L^p(K)$,  for every $k$ and every compact set $K\subset\R$.}
\end{remark}

\begin{lemma}\label{prop:dec} Let $\bar u:\R\mapsto\R$ be a  measurable function.
Then  $\ol u\in {\cal P}_{\alpha}$ if and only if it can be decomposed as (\ref{udec}), 
     where the $v_k$ satisfy (\ref{v00}) and (\ref{vk}).   
    The smallest constant $C$ for which (\ref{v00}) and (\ref{vk}) hold for every $k\geq 1$  is of the same order of $\Vert \ol u \Vert_{{\cal P}_\alpha}$.
\end{lemma}

{\bf Proof.}
{\bf 1.} Assume that $\ol u$ admits a decomposition as in (\ref{udec})-(\ref{vk}). 
We show that 
$$
\Vert  \ol u \Vert_{{\cal P}_{\alpha}} ~=~\O(1) \cdot C.
$$
Consider the case where $\lambda\in \,] 0,1]$ is of the form $\lambda = 2^{-q}$ for some integer $q\geq 1$. Then we set
$$
\wt v \,\doteq\,  \sum_{k =0}^q v_k
$$
and estimate 
$$
\TV\{ \wt v\}~ \leq ~\sum_{k = 0}^q \TV\{ v_k\}~ \leq~ C \cdot \sum_{k =0}^q 2^{(1-\alpha) k} ~= ~\O(1) \cdot C \cdot 2^{(1-\alpha) q}. 
$$
Moreover 
$$
\meas\bigl( \{ \ol u\ne \wt v \} \bigr)~ \leq ~C \cdot \sum_{k =q+1}^{+\infty} 2^{-\alpha k} ~=~ \O(1) \cdot C \cdot  2^{-\alpha q}.
$$
This proves
$$
\sup_{q \geq 1}~ d^{(2^{-q})}(\ol u, 0)~=~  \O(1) \cdot C.
$$
A simple argument now shows that the estimate 
holds when the supremum is taken over all $0 < \lambda \leq 1$. 

\v
{\bf 2.} Assume now that $\ol u \in {\cal P}_{\alpha}$. By the definition of $\Vert \cdot \Vert_{{\cal P}_{\alpha}}$ at (\ref{d66})-(\ref{npa}), choosing $\lambda = 2^{-k}$, for every $k \geq 0$ we obtain a function $u_k$ such that
      \begin{equation}\label{eq:ukdef}
   \TV\{u_k\}\, \leq \,\Vert \ol u \Vert_{{\cal P}_{\alpha}}\cdot 2^{(1-\alpha) k}, \qquad \meas \Big( \{\ol u\ne u_k\}\Big)~ \leq ~\Vert \ol u\Vert_{{\cal P}_{\alpha}}\cdot 
   2^{-\alpha k}.
    \end{equation}
  One can choose the functions $u_k$ so that they also satisfy
$$
\{\ol u\ne u_k\} \,\subseteq\, \{ \ol u\ne u_{k-1}\} \qquad \forall \; k \geq 1.
$$ 
    We define the functions $\{v_k\}_{k\geq 1}$ by setting
    \begin{equation}
      v_0 \doteq u_0, \qquad v_k \doteq u_{k}-u_{k-1} \qquad \forall \; k \geq 1.
    \end{equation}
  By  the second inequality in \eqref{eq:ukdef} it follows
    $$
    \lim_{k \to +\infty} u_k(x) = \ol u(x), 
    \qquad \text{pointwise for a.e. $x \in \mathbb R$}.
    $$
Using a telescopic sum one obtains
$$
\ol u(x) ~=~ \lim_{N \to +\infty} u_N(x) ~ =~ \lim_{N \to +\infty}  \sum_{k = 0}^{N} v_k(x) \qquad \text{pointwise for a.e. $x \in \mathbb R$}.
$$
We conclude by observing that
$$
\TV\{ v_k\}~ \leq~ \TV\{ u_k\} + \TV \{u_{k-1}\}~ \leq ~2\cdot \Vert \ol u \Vert_{{\cal P}_{\alpha}}\cdot 2^{(1-\alpha) k}
$$
and  
$$
\meas\big(\{u_k\ne u_{k-1}\} \big) ~\leq ~\meas \big(\{\ol u\ne u_{k-1}\} \big) ~  \leq~ 2^{\alpha} \cdot \Vert \ol u\Vert_{{\cal P}_{\alpha}} \cdot 2^{-\alpha k}.
$$
\endproof

We now show that one can choose the decomposition in Lemma \ref{prop:dec}
in such a way that all functions $v_k$ are one-sided Lipschitz.
\begin{lemma}\label{lemma:strong_dec} Consider any
function $\ol u\in{\cal P}_{\alpha}$.  Then, it is possible to choose the functions $v_k$ in (\ref{udec})--(\ref{vk}) in such a way that the additional
one-sided Lipschitz bound (\ref{oslip}) holds.

\end{lemma}

{\bf Proof.} 
{\bf 1.}
In the following, given $p>0$ and a function $f$, we denote by $\caL_p (f)$ its lower one-sided $p$-Lipschitz envelope:
$$
   \caL_p f (x)~\doteq ~{\rm sup} \Big\{ u(x) \,;~~ u:\mathbb R \to \mathbb R, \  u(y) \leq f(y)\,,~~u(y')-u(y)\leq p (y'-y) ~~\forall \; y,y'\in\R, ~~y <y' \Big\}.
    $$

\begin{figure}[ht]
  \centerline{\hbox{\includegraphics[width=8cm]{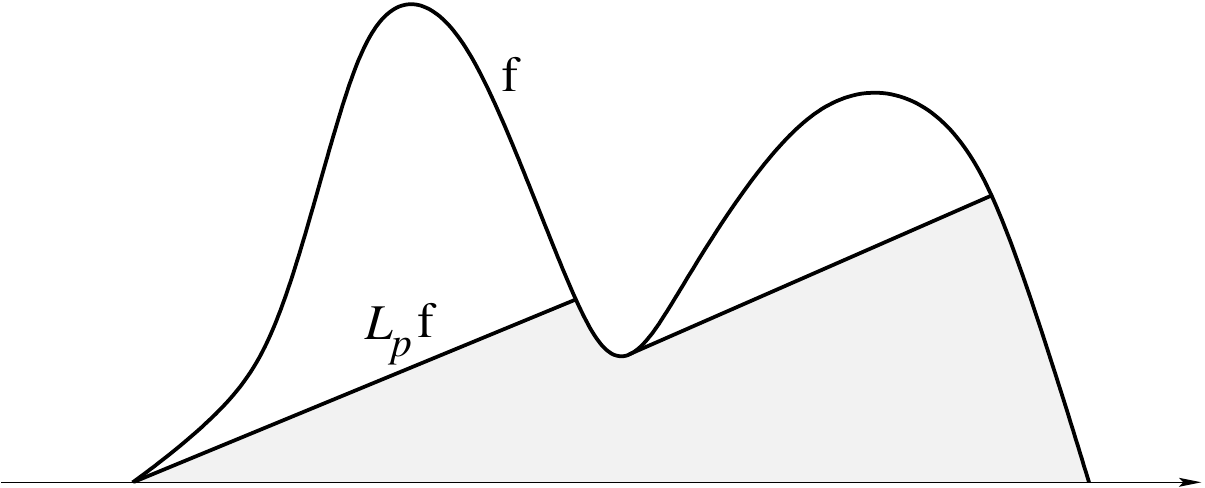}}}
  \caption{\small  The lower one-sided $p$-Lipschitz envelope $\caL_p f$.  Here the straight lines have slope $p$.  }\label{f:ag45}
\end{figure}


The function $\caL_p f$ is the largest one-sided $p$-Lipschitz function whose graph
lies below the graph of $f$. Denoting by $\TV^+\{g\}$ the positive variation of a function $g$
$$
\TV^+\{g\} \doteq \sup\! \sum_{x_0\le\cdots\le x_n}\big[g(x_{i})-g(x_{i-1})\big]^+, \qquad \big[g(x_{i})-g(x_{i-1})\big]^+ = \max\big\{0, \,g(x_{i})-g(x_{i-1})\big\},
$$
the following  relations  hold:
    \begin{equation}\label{eq:envelopeest}
   \left\{     \bega{l}
             \TV\!^+ \{\caL_p f\}~ \leq ~\TV\!^+ \{f\}, \\[2mm]
             \meas \Big\{ x \in \R~;~~\caL_p f(x) < f(x) \Big\} ~\leq~
             \frac{1}{p} \cdot \TV\!^+\{f\}, \\[2mm]
             \TV\!^+ \{f-\caL_p f\} ~\leq ~\TV\!^+ \{f\},
        \enda\right.
    \end{equation}
  Notice that the second inequality is a consequence of Riesz' sunrise lemma 
       (see for example \cite{KF}, p.319), while 
       the other two inequalities are straightforward.
\v
{\bf 2.} 
    Let $\{p_k^i\}_{i,k \in \mathbb N}$ be positive numbers (to be chosen later) such that 
    \begin{equation}\label{eq:pki}
    \lim_{i \to +\infty} p_{k}^i = +\infty \qquad \forall \; k \geq 0,
    \end{equation}
    and 
    consider the decomposition (\ref{udec}) constructed in   Lemma~\ref{prop:dec}.
   As an intermediate step,  we claim that for every $ k \geq 0$, the function $v_k$  
   can be further decomposed as 
   \bel{57}
    v_k(x) ~=~ \sum_{i=0}^{+\infty} v_k^i(x) \qquad \text{for a.e.}~x \in \R,
    \eeq
    where each $v_k^i$ is one-sided $p^i_k$-Lipschitz and satisfies 
\bel{58}
\begin{aligned}
\meas\bigl(\{v_k^i\ne 0 \}\bigr)~\leq ~C \cdot \frac{2^{k(1-\alpha)}}{p_k^{i-1}}, \qquad & \TV\!^+ \{v_k^i\}~ \leq~ C \cdot 2^{(1-\alpha) k}  \quad \text{if $i+k \neq 0$} ,\\
 &  \TV \{v_0^0\}~ \leq~ C,
\end{aligned}
\eeq
with  $C= \O(1)\cdot \Vert \ol u\Vert_{\P_{\alpha}}$, as in Lemma~\ref{prop:dec}.

For $k \geq 0$, we first prove the claim assuming $v_k\geq 0$.
Define 
    $$
     v^0_k \doteq \caL_{p_k^0} (v_k).
    $$
By \eqref{eq:envelopeest} one has
$$
\begin{aligned}
& \TV\!^+ \{  v^0_k\}~ \leq~ \TV\!^+ \{ v_k \}~\leq~  C \cdot 2^{(1-\alpha) k}, \qquad \text{for all $k \geq 0$}. \\
& \meas(\{v^0_k\ne 0\} )~ \leq \meas(\{ v_k\ne 0\} )~ \leq~C \cdot 2^{-\alpha k} \qquad \text{for all $k \geq 1$}.
\end{aligned}
$$

Setting $\psi_0 = v_k -v^0_k$ and using again \eqref{eq:envelopeest} we get
$$\meas( \{\psi_0\ne 0\})~ \leq ~ C \cdot \frac{1}{p_k^0}\cdot 2^{(1-\alpha) k}, \qquad \TV\!^+ \{ \psi_0\}~ \leq ~C \cdot 2^{(1-\alpha) k}.
$$
Defining $v_k^1 = \caL_{p_k^1} \psi_0$ we obtain 
$$\meas(\{ v_k^1\ne 0\})~ \leq~ C \cdot \frac{1}{p_k^0}\cdot 2^{(1-\alpha) k}, \qquad \TV\!^+ \{v_k^1\}~ \leq~ C \cdot 2^{(1-\alpha) k}.
$$
By induction, assume we are given $\psi_{i-2}$ and $v_k^{i-1} = \caL_{p_k^{i-1}} \psi_{i-2}$, $i \geq 2$, both with positive variation  $\leq C \cdot 2^{(1-\alpha) k}$. We then define  $\psi_{i-1} = \psi_{i-2}-v_k^{i-1}$ and $v_k^{i} = \caL_{p_k^i} \psi_{i-1}$. This yields
$$\meas(\{\psi_{i-1}\ne 0\} )~ \leq ~ C \cdot \frac{1}{p_k^{i-1}}\cdot 2^{(1-\alpha) k}, \qquad \TV\!^+ \{\psi_{i-1}\}~ \leq ~C \cdot 2^{(1-\alpha) k}
$$
and hence
\bel{Bi}
\meas(\{ v_k^i\ne 0\})~ \leq ~C \cdot \frac{2^{(1-\alpha)k}}{  p_k^{i-1}}, \qquad \TV\!^+ \{v_k^i\}~ \leq~ C \cdot 2^{(1-\alpha) k}.
\end{equation}
Therefore, by induction  \eqref{Bi} holds for every $i$ and every $k$. 

Finally, we prove that \eqref{57} holds. Indeed, by the second inequality in \eqref{eq:envelopeest} it follows
$$
{\rm meas} \, \left(\left\{v_k\ne \sum_{i=0}^{i^*}v_k^i \right\}\right)~ \leq ~\frac{1}{p_k^{i^*}} \cdot \TV \, v_k \qquad \forall \; i^* \geq 1.
$$
Letting $i^*\to +\infty$, this proves our claim in the case $v_k\geq 0$.
\v
{\bf 3.}
Next, we show how to handle the general case where $v_k = v_k^+ - v_k^-$
has a positive and a negative part.  We already know how to decompose the positive part $v_k^+$.

We treat the negative part $v_k^-$ in the same way, but using instead the lower one-sided Lipschitz envelope, defined by
$$
   \caL_p^- f (x)~\doteq ~{\rm sup} \Big\{ u(x) \,;~~ u(y) \leq f(y)\,,~~u(y')-u(y)\geq - p (y'-y) ~~\forall \; y,y'\in\R, ~~y <y' \Big\}.
    $$

%
This yields a decomposition 
$$
v_k^- = \sum_{i = 0}^{\infty} w_k^i
$$
where $w_k^i$ are positive one sided Lipschitz functions  which satisfy 
the same inequalities as in  (\ref{58}), and whose distributional derivatives satisfy $D w_i^k \geq -p_k^i$.  

Then
$$v_k(x)~=~v_k^+(x)-v_k^-(x)~=~\sum_{i = 0}^{\infty} v_k^i(x)+\sum_{i = 0}^{\infty} \bigl(-w_k^i(x)\bigr).$$
is the desired decomposition.

\v
{\bf 4.} We now conclude the proof of the lemma, relying on
(\ref{57})-(\ref{58}). 
Defining
\bel{59}
\wt v_q ~= ~\sum_{i= 0}^q v^{i}_{q-i}\,,
\eeq
we obtain
\bel{510}
    \ol u~=~ \sum_{q=0}^{+\infty} \wt v_q\,.
\eeq
Next, we choose 
$$p_k^i\, \doteq \,\frac{6}{\pi^2}\cdot \frac{ 2^{k+i}}{(i+2)^{2}}\,.$$
By (\ref{58}),  for every $q \geq 1$ we obtain
$$\meas\bigl(\{ \wt v_q \ne 0\}\bigr)~ \leq~\frac{\pi^2}{6} C \cdot 2^{-\alpha q} \sum_{i=0}^q \frac{2^q\cdot (i+1)^2}{2^{(1-\alpha) i}\cdot 2^{q-1}}~ \leq~ \frac{\pi^2}{3}\cdot C \cdot 2^{-\alpha q} \cdot \sum_{i=0}^q \frac{(i+1)^2}{2^{(1-\alpha) i}} ~= ~\O(1) \cdot C \cdot 2^{-\alpha q}.
$$
Moreover, always for $q \geq 1$, the one-sided Lipschitz constant for $\Tilde v_q$ is estimated by
$$
D \wt v_q ~\leq~ \frac{6}{\pi^2}\sum_{i=0}^q p_{q-i}^i ~\leq ~2^q \frac{6}{\pi^2}\sum_{i=0}^q \frac{1}{(i+1)^2}~ \leq~   2^q.
$$  
The one-sided Lipschitz property and the estimate on the support of $\widetilde v_q$ readily imply $$\TV \{\widetilde v_q \}~=~ \O(1)\cdot C \cdot 2^{(1-\alpha)q}.$$ If $q = 0$, by definition we have $\widetilde v_q = v_0^0$, so that by \eqref{58}
$$
\TV \,\{\widetilde v_0\} \, \leq \, C, \qquad D \widetilde v_0\, \leq \, p_0^0\,  \leq\,  1 .
$$
The conclusion of the lemma is achieved by
renaming $v_k \doteq \widetilde v_q$, with $q = k$.
\endproof

\begin{lemma}\label{l:54}
 For every $k \geq 0$, the function $ v_k$ constructed in Lemma \ref{lemma:strong_dec}
 satisfies 3 of Theorem \ref{thm:dec}.
    \end{lemma}

{\bf Proof.}
    Since $v_k$ is one-sided Lipschitz, the set $\{v_k\ne 0\}$ has at most countably many connected components,  that we denote by $\{I_k^p\}_{p \geq 1}$.
Consider then the restrictions $v_k^p \doteq  v_k\Big|_{I_k^p}$. On every interval $I_k^p$,  having length $\ell_k^p$, since $v_k$ is one-sided $ 2^k$-Lipschitz, one has 
$$
| v_k(x)| ~\leq ~ 2^k\ell_k^p~ =~ h_k^p, \qquad \forall \; x \in I_k^p
$$
Therefore
$$
\sum_{p \geq 1} \ell_k^p ~=~ \meas(\{v_k\ne 0\})~ \leq ~  C 2^{-\alpha k},
$$
which proves \eqref{eq:rectsum}.
    \endproof

\section{Decay rate of the Total Variation} 
\label{s:6}

In this section we prove that if $1/2< \alpha <1$, then the conjectured decay of the total variation with rate $t^{\alpha-1}$ holds.

\begin{theorem}\label{t:61} Consider a bounded, compactly supported initial datum 
  $\ol u \in {\cal P}_{\alpha}$, with $1/2 <\alpha <1$. Then the solution to Burgers' equation (\ref{Bur}) satisfies
\bel{adecay}
    \limsup_{t \to 0+} \;\Big(t^{1-\alpha}\, \cdot  TV\, \bigl\{ S_t \ol u \bigr\}\Big)\; \leq\;  C_0\,\frac{ \Vert \ol u\Vert_{{\cal P}_{\alpha}}}{2\alpha-1}\,,
 \eeq
     where $C_0$ is some absolute constant. 
\end{theorem}

%
%
{\bf Proof.} {\bf 1.} 
Let $u=u(t,x)$ be a solution of Burgers' equation and let  $t>0$ be given. Denote by $J_t\subset \R$ the jump set of $u(t,\cdot)$. By the Lax-Oleinik formula, for every $x\in \R\setminus J_t$ there is a unique backward characteristic from the point $(t,x)$ along which the solution is constant: namely 
\bel{charBurgers}
u(t,x)=u(s, x-u(t,x)(t-s)) \qquad \forall \, s \in \,]0,t].
\eeq
If \eqref{charBurgers} holds we say that the couple 
$$(x_0,v)~\doteq~ \bigl(x-u(t,x)t\,, ~u(t,x)\bigr)$$ {\it survives up to time $t$}, and we denote by $\mathcal Q{(t)}$ the set of couples which survive up to time $t$. 
The point $x_0$ is the starting point of the characteristic passing through $(t,x)$ and $v$ is the value of $u$ along the characteristic. Notice however, for example in the case of a centered rarefaction, that we can have $\ol u(x_0) \ne v$. Indeed, \eqref{charBurgers} does not extend to $t=0$, in general.
We will estimate $\TV\{S_t \ol u\}$ by means of the equality
\bel{traceback}
\TV\bigl\{S_t \ol u\,;~ \R\bigr\} ~= ~\TV\bigl\{S_t \ol u\,;~ \R \setminus J_t\bigr\}~=~ \sup \sum_{\substack{x_1\le\cdots\le x_n,\\  (x_i,v_i)\in \mathcal{Q}(t)}}|v_i-v_{i-1}|,
\eeq
where in the last sum we assume that if $x_{i-1}=x_i$, then $v_{i-1}\le v_i$.
By the Lax formula, the constraint $(x_0,u_0)\in \mathcal Q(t)$ is satisfied if and only if
\begin{equation}
    \begin{aligned}
        \int_{x_0}^y \ol u(z) -\Big[u_0 -\frac{1}{t}(z-x_0) \Big]\,d z\, \geq\, 0 \qquad \forall \; y \geq x_0\,,\\
        \int_{y}^{x_0} \ol u(z) -\Big[u_0 -\frac{1}{t}(z-x_0) \Big] \,d z\, \leq\, 0 \qquad \forall \; y \leq x_0\,.
    \end{aligned}
\end{equation}
Equivalently:
\begin{equation}\label{eq:rightsurv}
 \int_{x_0}^y \Big[\ol u(z)-\Big(u_0 -\frac{z-x_0}{t}\Big)\Big]^+ \,d z ~
 \geq~  \int_{x_0}^y \Big[\ol u(z)-\Big(u_0 -\frac{z-x_0}{t}\Big)\Big]^- \,d z   
 \qquad \forall \; y \geq x_0\,,
        \end{equation}
    \begin{equation}\label{eq:leftsurv}
         \int_{y}^{x_0} \Big[\ol u(z)-\Big(u_0 -\frac{z-x_0}{t}\Big)\Big]^+ \,d z ~\leq ~ \int_{y}^{x_0} \Big[\ol u(z)-\Big(u_0 -\frac{z-x_0}{t}\Big)\Big]^- \,d z  \qquad \forall \; y \leq x_0\,,
\end{equation}
where we used the notation
$$
 [z]^+ \doteq \max \{z, 0\}, \qquad [z]^- \doteq -\min \{z, 0\}.
$$
The interpretation of \eqref{eq:rightsurv} is that for every $y > x_0$ the area of the hypograph of $\ol u$ in $[x_0,y]$ that lies above the line passing through $(x_0, u_0)$ with slope $-1/t$ must be bigger then the area of the epigraph lying below the same line. 
Analogously, the interpretation of \eqref{eq:leftsurv} is that for every $y < x_0$ the area of the hypograph of $\ol u$ in $[x_0,y]$ that lies above the line passing through $(x,u_0)$ with slope $-1/t$ must be smaller then the area of the epigraph lying below the same line.

\v
{\bf 2.} 
It suffices to prove the decay estimate (\ref{adecay}) 
    for all times of the form $t = 2^{-k}$, $k\geq 1$.
 We thus need to estimate the quantity
    $$
        \limsup_{k \to +\infty} \; 2^{(\alpha-1) k}\, \TV \bigl\{  S_{2^{-k}} \ol u
        \bigr\} 
    $$
    In the following we fix a time $t = 2^{-k}$ and show that 
\begin{equation}\label{eq:t2kthm}
\TV \bigl\{  S_{2^{-k}} \ol u
        \bigr\} ~\leq~ c \cdot  {2^{(1-\alpha)k}\over 2\alpha-1}\cdot  \Vert \ol u\Vert_{{\cal P}_{\alpha}}\,,
 \end{equation}
where $c$ is some constant depending only on $\Vert \ol u\Vert_{\mathbf L^{\infty}}$.

Let $\ol u = \sum_{q=0}^\infty v_q$ be a decomposition satisfying all the properties
listed in Theorem~\ref{thm:dec}.
We write $\ol u$ as the sum of two terms:
\bel{ovud}
\ol u ~= ~\sum_{q=0}^{k-1} v_q + \sum_{q = k}^{\infty} v_q ~\doteq ~\widetilde u_k + \widehat u_k.
\eeq
We regard the function $\wt u_k$ as the  {\it regular part of $\ol u$},  in the sense that it has total variation that is bounded, of size $\O(1)\cdot 2^{(1-\alpha) k}$.  In fact,  since each $v_q$ is one-sided Lipschitz with constant $2^q$ and with total variation bounded by $C \,2^{(1-\alpha) q}$,  the function $\wt u_k$ is one-sided Lipschitz with constant $2^k$ and with total variation bounded by $C 2^{(1-\alpha) k}$.   

We recall that $C$ is the constant coming from Lemmas \ref{prop:dec} and \ref{lemma:strong_dec}, of the same order of $\Vert \ol u \Vert_{\P_{\alpha}}$.
\v
{\bf 3.} 
To simplify the exposition, we first give a proof under two additional assumptions:
\begin{enumerate}
    \item[(H1)] the regular part $\wt u_k$ is zero, i.e.
    \begin{equation}\label{eq:simpl_proof}
\ol u\,= \,\sum_{q = k}^{+\infty} v_q \,=\, \widehat u_k.
\end{equation}
    \item [(H2)] All functions $v_q$ are positive:
    $$
    v_q(x)\,\geq\, 0 \qquad \forall \; x \in \mathbb R.
    $$
    This implies that all the functions $v_q^p$ constructed in Lemma \ref{l:54} are positive as well,
\end{enumerate}

Assuming (H1) and (H2),
let $x_0 \leq x_1 \leq \cdots \leq x_n$ and $v_0,\ldots, v_n$,  be such that $(x_i, v_i) \in \Q({2^{-k}})$. 
According to \eqref{traceback}, it suffices to estimate the total variation over these points:
$$
\sum_{i=1}^n |v_i-v_{i-1}|.
$$
Since $\ol u$ is compactly supported, it is enough to estimate the negative variation, namely
\bel{negvar}
\sum_{i=1}^n \bigl[v_i-v_{i-1}\bigr]^-.
\eeq
The set of downward jumps
$$
{\cal N} ~\doteq ~\Big\{ i \in \{1, \ldots, n\} \,;~~ 
v_{i}\,<\, v_{i-1}\Big\}.
$$
can be partitioned as ${\cal N} = {\cal I} \cup {\cal J}$, where
$$
{\cal I} \doteq \Big\{ i \in {\cal N} \,;~~]x_{i-1}, x_i[ \,\subseteq \{ \widehat u_k \ne 0 \}\Big\}, \qquad {\cal J} \doteq {\cal N} \setminus {\cal I}.
$$
In the next two steps, the negative variation (\ref{negvar}) will be estimated 
by considering the  terms $i\in {\cal I}$ and $i\in {\cal J}$ separately.
\v
{\bf 4.} Let $i\in {\cal I}$. 
Since the characteristics starting at $x_i, x_{i+1}$ 
do not cross up to time  $t=2^{-k}$, this implies
$$
v_{i-1}-v_{i} ~\leq ~(x_{i}-x_{i-1}) \cdot 2^k.
$$
Summing over ${\cal I}$ one obtains
\bel{decay1}
\sum_{i \in {\cal I}} \bigl[v_{i-1}-v_{i}
\bigr]
~\leq ~2^k \,\sum_{i \in {\cal I}} (x_i-x_{i-1}) \leq 2^k \cdot \meas\bigl( \{ \widehat u_k \ne 0\}\bigr)~ \leq ~2C \cdot 2^k \cdot 2^{-\alpha k} ~= ~c_1 \,C \cdot 2^{(1-\alpha) k},
\eeq
where we used the inequalities
$$
\meas\big(\{\widehat u_k\ne 0\}\big) ~\leq ~\sum_{q= k}^{\infty}\meas\big( \{v_q\ne 0\}\big) ~\leq ~C \sum_{q= k}^{\infty} 2^{-\alpha q}~ \leq ~\frac{2^{\alpha}}{2^{\alpha}-1}C \cdot 2^{-\alpha k}.
$$

\begin{figure}
    \centering
\tikzset{every picture/.style={line width=0.75pt}} 

\begin{tikzpicture}[x=0.75pt,y=0.75pt,yscale=-1.4,xscale=1.4]

\draw    (90.33,200.33) -- (345.67,201.67) ;
\draw  [color={rgb, 255:red, 155; green, 155; blue, 155 }  ,draw opacity=1 ][line width=0.75] [line join = round][line cap = round] (97.67,199) .. controls (97.67,198.01) and (99.63,198.37) .. (100.33,197.67) .. controls (101.9,196.1) and (103.3,194.3) .. (104.33,192.33) .. controls (108.31,184.76) and (111.65,176.87) .. (115,169) .. controls (125.46,144.4) and (128.18,123.01) .. (134.33,95.67) .. controls (135.91,88.67) and (138.83,70.6) .. (146.33,69.67) .. controls (148.85,69.35) and (150.09,84.74) .. (150.33,87) .. controls (153.55,116.87) and (156.8,146.8) .. (162.33,176.33) .. controls (163.17,180.81) and (167.2,201) .. (173.67,201) ;
\draw  [color={rgb, 255:red, 155; green, 155; blue, 155 }  ,draw opacity=1 ][line width=0.75] [line join = round][line cap = round] (187,199.67) .. controls (193.02,185.61) and (192.92,171.11) .. (196.33,156.33) .. controls (198.53,146.81) and (196.26,132.7) .. (205,128.33) .. controls (206.65,127.51) and (208.06,130.04) .. (208.33,131) .. controls (210.54,138.74) and (210.96,147.77) .. (212.33,155.67) .. controls (214.42,167.68) and (217.01,180.9) .. (221.67,192.33) .. controls (222.72,194.91) and (225.67,196.88) .. (225.67,199.67) ;
\draw [color={rgb, 255:red, 155; green, 155; blue, 155 }  ,draw opacity=1 ][line width=0.75] [line join = round][line cap = round]   (234.33,200.33) .. controls (237.6,200.33) and (243.91,174.01) .. (245,169) .. controls (246.6,161.65) and (246.6,152.98) .. (247.67,145.67) .. controls (250.55,125.99) and (253.92,106.53) .. (257,87) .. controls (257.86,81.52) and (260.62,70.87) .. (267,71.67) .. controls (269.4,71.97) and (268.95,81.11) .. (269,81.67) .. controls (269.77,90.47) and (272.44,110.32) .. (273,115.67) .. controls (274.82,132.99) and (277.63,169.91) .. (288.33,183.67) .. controls (293.02,189.69) and (295.58,178.91) .. (297.67,181) .. controls (301.69,185.02) and (300.91,201) .. (307,201) ;
\draw  [dash pattern={on 4.5pt off 4.5pt}]  (74.33,93) -- (316.33,218.33) ;
\draw  [dash pattern={on 0.84pt off 2.51pt}]  (86.67,120) -- (341,121) ;
\draw  [dash pattern={on 0.84pt off 2.51pt}]  (87.1,173.05) -- (339.67,171.67) ;
\draw  [fill={rgb, 255:red, 245; green, 166; blue, 35 }  ,fill opacity=0.64 ] (228.54,172.44) -- (129.63,173.25) -- (129.2,121.01) -- cycle ;
\draw    (146.67,183.67) -- (182.63,243.29) ;
\draw [shift={(183.67,245)}, rotate = 238.9] [color={rgb, 255:red, 0; green, 0; blue, 0 }  ][line width=0.75]    (10.93,-3.29) .. controls (6.95,-1.4) and (3.31,-0.3) .. (0,0) .. controls (3.31,0.3) and (6.95,1.4) .. (10.93,3.29)   ;
\draw    (212.67,191.67) -- (204.65,242.36) ;
\draw [shift={(204.33,244.33)}, rotate = 278.99] [color={rgb, 255:red, 0; green, 0; blue, 0 }  ][line width=0.75]    (10.93,-3.29) .. controls (6.95,-1.4) and (3.31,-0.3) .. (0,0) .. controls (3.31,0.3) and (6.95,1.4) .. (10.93,3.29)   ;
\draw    (91.33,92) -- (103.06,54.91) ;
\draw [shift={(103.67,53)}, rotate = 107.55] [color={rgb, 255:red, 0; green, 0; blue, 0 }  ][line width=0.75]    (10.93,-3.29) .. controls (6.95,-1.4) and (3.31,-0.3) .. (0,0) .. controls (3.31,0.3) and (6.95,1.4) .. (10.93,3.29)   ;
\draw  [dash pattern={on 0.84pt off 2.51pt}]  (129.2,121.01) -- (129.67,200.33) ;
\draw  [dash pattern={on 0.84pt off 2.51pt}]  (245,173) -- (245,201) ;
\draw    (77.48,125.5) -- (77.19,167.83) ;
\draw [shift={(77.17,170.83)}, rotate = 270.4] [fill={rgb, 255:red, 0; green, 0; blue, 0 }  ][line width=0.08]  [draw opacity=0] (5.36,-2.57) -- (0,0) -- (5.36,2.57) -- cycle    ;
\draw [shift={(77.5,122.5)}, rotate = 90.4] [fill={rgb, 255:red, 0; green, 0; blue, 0 }  ][line width=0.08]  [draw opacity=0] (5.36,-2.57) -- (0,0) -- (5.36,2.57) -- cycle    ;

\draw (162.7,151.9) node [anchor=north west][inner sep=0.75pt]  
{$T_{i}$};
\draw (90.2,36.4) node [anchor=north west][inner sep=0.75pt]  
{$slope\ = \ -2^{k}$};
\draw (173.2,250.4) node [anchor=north west][inner sep=0.75pt]  
{$bumps\ v_{q}^{p}$};
\draw (120,207.3) node [anchor=north west][inner sep=0.75pt]  
{$x_{i}{}_{-1}$};
\draw (240,209.3) node [anchor=north west][inner sep=0.75pt]  
{$x_{i}$};
\draw (61,141.8) node [anchor=north west][inner sep=0.75pt] 
 {$\delta _{i}$};
\draw (316.5,103.3) node [anchor=north west][inner sep=0.75pt]  
{$v_{i-1}$};
\draw (317,154.8) node [anchor=north west][inner sep=0.75pt]  
{$v_i$};
\draw (155.33,75.4) node [anchor=north west][inner sep=0.75pt]  
 {$\widehat{u}_{k}$};
\end{tikzpicture}
    \caption{\small The configuration considered in the estimate for $i \in {\cal J}$, in 
    step {\bf 5} of the proof. }
    \label{fig:triangle}
\end{figure}
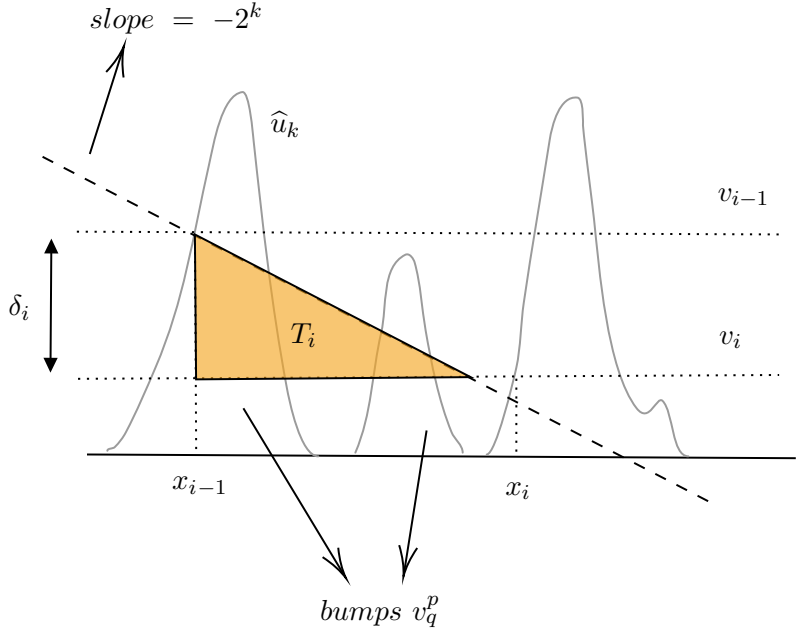

\v {\bf 5.} Next, consider the case where $i\in {\cal J}$. 
As shown in  Fig.~\ref{fig:triangle}, consider the triangle 
$$
\T_i ~\doteq ~\Big\{ (x, v) \in \mathbb R^2 \,;~~x_{i-1}<x< x_i, \quad v_{i}\leq v \leq v_{i-1}-2^k  (x-x_{i-1}) \Big\}.
$$
with height  $\delta_i \doteq v_{i-1} - v_{i}$ 
and base of length $\delta_i2^{-k}$.

In the following, we denote by $\widehat U_k \subset \mathbb R^2$  the region below the
graph  of $\widehat u_k$:
$$
\widehat U_k ~\doteq~ \Big\{ (x, v) \in \mathbb R^2 \,;~~0 \leq v \leq \widehat u_k(x) \Big\}~ \subset~ \mathbb R^2.
$$
By \eqref{eq:rightsurv}, the fact that $(x_{i-1},v_{i-1}) \in \Q(2^{-k})$,  implies that the area of the triangle $\T_i$ is bounded by the area of $\widehat U$ in the strip $[x_{i-1}, x_i]\times \mathbb R$: 
$$
\meas(\T_i) ~=~2^{-(k+1)}\delta_i^2~\leq ~\meas \Big( \widehat U_k \cap \bigl( [x_{i-1}, x_i]\times \mathbb R 
\bigr)\Big).
$$
 This implies 
\begin{equation}\label{eq:triest}
\delta_i ~\leq~ 2^{\frac{k+1}{2}}\cdot \meas 
\Big( \widehat U_k \cap \bigl( [x_{i-1}, x_i]\times \mathbb R 
\bigr)\Big)^{1/2}.
\end{equation}
For each $i \in {\cal J}$, we now consider the set of indices ${\cal Z}(i) \subset \mathbb N \times \mathbb N$ defined by 
$$
{\cal Z}(i) \,\doteq \,\Big\{ (p, q) \in \mathbb N^2 \,;~~ q \geq k, \quad {\rm supp}\, v_q^p \,\cap \,]x_{i-1}, x_i[\,  \neq \emptyset\Big\}.
$$
This is the set of all functions $v_q^p$ in the decomposition of  $\widehat u_k$ whose support intersects the open  interval $\,]x_{i-1},x_i[\,$. 

At this stage, we make an important observation: 
\begi
\item {\it For any couple $(p,q)$
with $q \geq k$, there can be at most two indices $i \in {\cal J}$ 
such that $(p, q) \in {\cal Z}(i)$}
\endi
Indeed, assume that this were not the case, i.e. for some $i_1< i_2 <i_3$ one had $(p, q) \in {\cal Z}(i_1) \cap {\cal Z}(i_2) \cap {\cal Z}(i_3)$.
Since the set $\{v_q^p\ne 0\}$ is connected, we would have
$$
]x_{i_2-1}, x_{i_2}[~ \subset ~\{v_q^p\ne 0\}~ \subseteq~\{\widehat  u_k\ne 0\}.
$$
But this is a contradiction because $i \notin {\cal I}$.

As in Lemma~\ref{l:54}, call $\ell_q^p \doteq \meas\bigl({\rm supp}\, v_q^p\bigr)$.
By \eqref{eq:triest} it follows
$$
\delta_i ~ \leq ~\sqrt{2} \cdot 2^{k/2}\cdot \Bigg(  \,\sum_{(p,q) \in {\cal Z}(i)} h_q^p \cdot \ell_q^p\, \Bigg)^{1/2} \leq  ~ \sqrt{2}\ \cdot 2^{k/2}\cdot  \sum_{(p,q) \in {\cal Z}(i)} \Big( h_q^p \cdot \ell_q^p\, \Big)^{1/2} .
$$
Summing over $i\in {\cal J}$, and using the fact that each couple $(p,q)$ 
can appear in the sum at most twice,  we obtain 
\begin{equation}\label{eq:twoiineq}
\begin{aligned}
    \sum_{i \in {\cal J}} \delta_i & \leq ~  \sqrt{2}\ \cdot 2^{k/2}\cdot  \sum_{i \in {\cal J}} \, \sum_{(p,q) \in {\cal Z}(i)} \Big( h_q^p \cdot \ell_q^p\, \Big)^{1/2} \\
    & \leq  ~2^{3/2}\ \cdot 2^{k/2}  \cdot \sum_{q = k}^{\infty}\, \sum_{p\in \mathbb N}\Big( h_q^p \cdot \ell_q^p\, \Big)^{1/2}.
    \end{aligned}
\end{equation}
Observing that
$$
\Big(\ell_q^p\cdot  h_q^p\Big)^{1/2} =~ 2^{q/2}\ell_q^p
$$
and using \eqref{eq:rectsum},  from   \eqref{eq:twoiineq} we obtain
\bel{decay2} 
\bega{rl}\ds
\sum_{i \in {\cal J}} \delta_i &\ds \leq ~   2^{3/2} \cdot 2^{k/2} \cdot \sum_{q=  k}^{\infty} \, \sum_{p \in \mathbb N}  \, 2^{q/2}\ell_q^p ~
~ \ds\leq~ 2^{3/2} \cdot  C \cdot  2^{k/2} \cdot \sum_{q =k}^{\infty}\,  2^{q/2}  2^{-\alpha q}\\[4mm]
 & \ds\leq  ~ c_2 \,C \cdot \frac{1}{1-2\alpha} 2^{k/2} \cdot 2^{(1/2-\alpha) k}~ =~ c_2 \,C \cdot \frac{1}{2\alpha-1}\cdot 2^{(1-\alpha) k},
\enda
\eeq
where $c_2$ is another absolute constant.
Combining (\ref{decay1}) with (\ref{decay2}), we obtain the desired decay rate,
under the additional assumptions (H1)-(H2).
\v
{\bf 6.} In the remaining steps we complete the proof of the theorem, removing the assumptions 
(H1)-(H2).   

Recalling the decomposition (\ref{ovud}), we observe that
the function $\wt u_k$ is one-sided $2^{k}$-Lipschitz, because each $v_q$ is 
one-sided $2^q$ Lipschitz.

Let $x_0 \leq x_1 \leq \ldots \leq x_n$ and $v_i$, $i=0, \ldots,  n$ be such that $(x_i, v_i) \in \Q({2^{-k}})$. 
As before, it suffices to estimate the negative variation, i.e.
$$
\sum_{i \in {\cal N}} \bigl[v_{i-1}-v_{i}\bigr],
$$
where 
$$
{\cal N}\, \doteq\, \Big\{ i \in \{1, \ldots, n\} \,;~~v_{i} < v_{i-1}\Big\}.
$$
We partition the above set of indices as ${\cal N} = {\cal I} \cup {\cal J}$, where
$$
{\cal I} \doteq \Bigg\{ i \in {\cal N} \,;~~]x_{i-1}, x_i[~ \subset\, \bigcup_{q\geq k} \{v_q \ne 0\} \Bigg\}, \qquad {\cal J} \doteq {\cal N} \setminus {\cal I}.
$$
Set 
$$
\delta_i \doteq v_{i-1}-v_{i} \qquad \forall \; i \in {\cal N}.
$$
We further partition  ${\cal J} = {\cal J}_1 \cup {\cal J}_2$, by setting
$$
     {\cal J}_1 \doteq \Big\{ i \in {\cal J} \; \mid \; \wt u_k(x_{i-1})\geq v_{i-1}- \frac{\delta_i}{3} \quad \text{and} \quad \wt u_k(x_i) \leq v_{i}+\frac{\delta_i}{3} \Big\},\qquad 
{\cal J}_2 \doteq {\cal J} \setminus {\cal J}_1\, .
$$
We will estimate the quantity 
$$
\sum_{i \in {\cal N}} \delta_i~ = ~\sum_{i \in {\cal I}} \delta_i+ \sum_{i \in {\cal J}_1} \delta_i + \sum_{i \in {\cal J}_2} \delta_i
$$
by providing a bound on each term on the right hand side,
 in the following three steps.

\v
{\bf 7. (Estimate of the sum over ${\cal I}$)}. With exactly the same argument as in Step 4., we obtain the estimate
\begin{equation}\label{es61}
 \sum_{i \in {\cal I}} \bigl[v_{i-1}-v_{i} \bigr]~\leq c_1 \cdot C  \cdot 2^{(1-\alpha) k}
\end{equation}
where $c_1$ is an absolute constant.

\begin{figure}
    \centering
    \tikzset{every picture/.style={line width=0.75pt}} 

\begin{tikzpicture}[x=0.7pt,y=0.7pt,yscale=-1,xscale=1]

\draw    (140.33,201.67) -- (140.33,49.67) ;
\draw    (441,50.33) -- (440.33,200.33) ;
\draw  [dash pattern={on 0.84pt off 2.51pt}]  (142.67,99.67) -- (439,99.67) ;
\draw  [dash pattern={on 0.84pt off 2.51pt}]  (142.67,150) -- (439,150.33) ;
\draw    (478.32,55.33) -- (477.68,194.67) ;
\draw [shift={(477.67,197.67)}, rotate = 270.26] [fill={rgb, 255:red, 0; green, 0; blue, 0 }  ][line width=0.08]  [draw opacity=0] (5.36,-2.57) -- (0,0) -- (5.36,2.57) -- cycle    ;
\draw [shift={(478.33,52.33)}, rotate = 90.26] [fill={rgb, 255:red, 0; green, 0; blue, 0 }  ][line width=0.08]  [draw opacity=0] (5.36,-2.57) -- (0,0) -- (5.36,2.57) -- cycle    ;
\draw  [line width=0.75] [line join = round][line cap = round] (124.33,88.33) .. controls (135.79,78.78) and (154.9,74.52) .. (169,71) .. controls (198.5,63.62) and (238.32,74.39) .. (263,91.67) .. controls (288.01,109.18) and (309.25,140.63) .. (333.67,160.33) .. controls (355.79,178.19) and (380.55,185.12) .. (408.33,187) .. controls (425.83,188.18) and (440.51,187.83) .. (455.67,179.67) ;

\draw (489.33,117.07) node [anchor=north west][inner sep=0.75pt]    {$\delta _{i}$};
\draw (114.33,53.73) node [anchor=north west][inner sep=0.75pt] 
 {$\frac{1}{3} \delta _{i}$};
\draw (113,158.07) node [anchor=north west][inner sep=0.75pt]  
{$\frac{1}{3} \delta _{i}$};
\draw (116.67,108.4) node [anchor=north west][inner sep=0.75pt]  
{$\frac{1}{3} \delta _{i}~$};
\draw (218.67,41.4) node [anchor=north west][inner sep=0.75pt]    {$\widetilde{u}_{k}$};
\end{tikzpicture}
    \caption{\small Illustration of the case $i \in {\cal J}_1$. At least one third of the variation is due to the regular part $\wt u_k$.}
    \label{fig:delta3}
\end{figure}
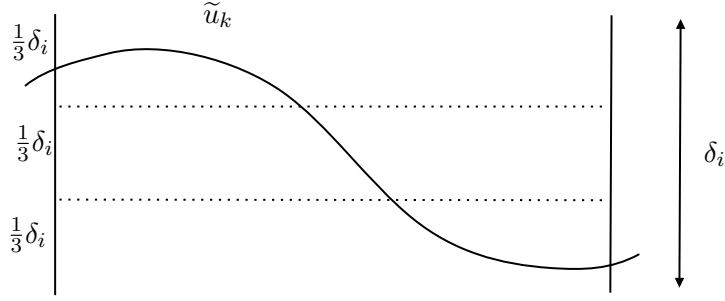

\v
{\bf 8. (Estimate of the sum over ${\cal J}_1$).}  In this case 
(see Figure~\ref{fig:delta3}),  by definition of ${\cal J}_1$ the variation of 
$\wt u_k$ on $[x_{i-1}, x_i]$  is at least one third of $\delta_i$. 
Therefore the jump $\delta_i$ is controlled by the variation of $\wt u_k$, which 
is the regular part. More precisely, from the definition of ${\cal J}_1$ it follows
$$
\begin{aligned}
   \delta_i ~\doteq ~ v_{i-1} - v_{i} ~\leq ~\Big(\wt u_k(x_{i-1}) + \frac{\delta_i}{3}\Big) -\Big(\wt u_k(x_i)-\frac{\delta_i}{3}\Big) ~ \leq ~\wt u_k(x_{i-1}) - \wt u_k(x_i)+ \frac{2\delta_i}{3}\,,
\end{aligned}
$$
and therefore
 $$
 \delta_i ~\leq~ 3 \cdot \big( \wt u_k(x_{i-1}) - \wt u_k(x_i)\big).
 $$
Summing over $i \in {\cal J}_1$ we obtain
\bel{es62}
    \sum_{i \in {\cal J}_1} \delta_i ~\leq ~3\sum_{i \in {\cal J}_1} \, \big( \wt u_k(x_{i-1}) - \wt u_k(x_i)\big) ~ \leq~ 3\,\TV \{ \wt u_k\}~ \leq ~c_2\cdot C \cdot 2^{(1-\alpha) k}.
\end{equation}

\v

{\bf 9. (Estimate of the sum over ${\cal J}_2$).} 
The idea here is that we reduced to a situation where we can proceed as in the simplified case of the previous section, up to a modification of the definition of the triangles $\T_i$ 
that takes into account the presence of $\wt u_k$, which is one-sided $2^k$-Lipschitz.
Since $i \in {\cal J}_2$, at least one of the following inequalities is true:
\bel{in12}
 \wt u_k(x_{i-1})<  v_{i-1}- \frac{\delta_i}{3} \quad\qquad \text{or} 
 \quad \qquad \wt u_k(x_i) > v_{i} +\frac{\delta_i}{3}.
\eeq
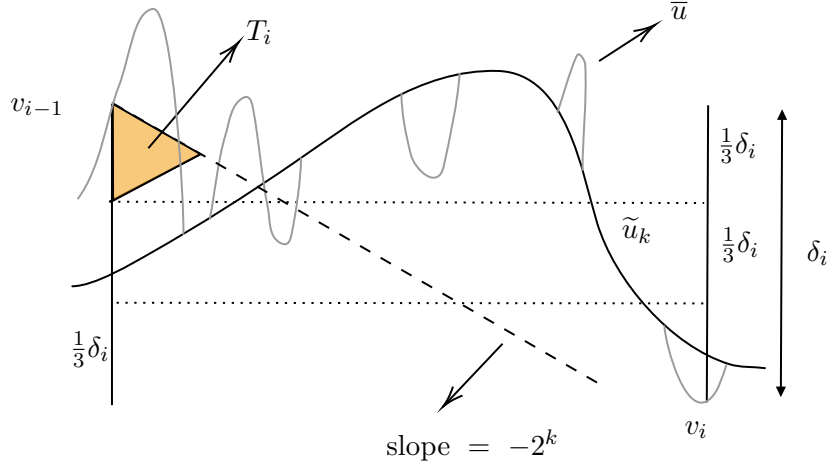
\begin{figure}
    \centering
    \tikzset{every picture/.style={line width=0.75pt}} 

\begin{tikzpicture}[x=0.75pt,y=0.75pt,yscale=-1,xscale=1]

\draw    (160.33,221.67) -- (160.33,69.67) ;
\draw    (461,70.33) -- (460.33,220.33) ;
\draw  [dash pattern={on 0.84pt off 2.51pt}]  (159,119) -- (459,119.67) ;
\draw  [dash pattern={on 0.84pt off 2.51pt}]  (162.67,170) -- (459,170.33) ;
\draw    (498.32,75.33) -- (497.68,214.67) ;
\draw [shift={(497.67,217.67)}, rotate = 270.26] [fill={rgb, 255:red, 0; green, 0; blue, 0 }  ][line width=0.08]  [draw opacity=0] (5.36,-2.57) -- (0,0) -- (5.36,2.57) -- cycle    ;
\draw [shift={(498.33,72.33)}, rotate = 90.26] [fill={rgb, 255:red, 0; green, 0; blue, 0 }  ][line width=0.08]  [draw opacity=0] (5.36,-2.57) -- (0,0) -- (5.36,2.57) -- cycle    ;
\draw  [line width=0.75] [line join = round][line cap = round] (140.33,161.67) .. controls (151.45,161.67) and (170.74,150.07) .. (179.67,145) .. controls (204.34,130.99) and (228.27,115.68) .. (251.67,99.67) .. controls (282.71,78.43) and (315.28,51.87) .. (355,53) .. controls (389.93,54) and (397.54,105.16) .. (405,129.67) .. controls (413.88,158.86) and (442.03,191.8) .. (471.67,201) .. controls (477.43,202.79) and (483.74,202.26) .. (489.67,203) ;
\draw  [color={rgb, 255:red, 155; green, 155; blue, 155 }  ,draw opacity=1 ][line width=0.75] [line join = round][line cap = round] (406.33,169.67) .. controls (406.33,169.67) and (406.33,169.67) .. (406.33,169.67) ;
\draw  [color={rgb, 255:red, 155; green, 155; blue, 155 }  ,draw opacity=1 ][line width=0.75] [line join = round][line cap = round] (209.67,127.67) .. controls (212.45,113.74) and (213.54,98.43) .. (217.67,85) .. controls (219.53,78.94) and (220.31,69.2) .. (227,66.33) .. controls (230.66,64.76) and (231.37,71.35) .. (231.67,73) .. controls (234.73,90.04) and (235.7,105.84) .. (239,122.33) .. controls (240.11,127.87) and (239.14,136.41) .. (245,139.67) .. controls (248.98,141.88) and (250.11,139.32) .. (251,135) .. controls (253.08,124.88) and (255.67,107.27) .. (255.67,97) ;
\draw  [color={rgb, 255:red, 155; green, 155; blue, 155 }  ,draw opacity=1 ][line width=0.75] [line join = round][line cap = round] (306.33,64.33) .. controls (306.33,74.67) and (308.54,84.96) .. (311,95) .. controls (312.97,103.06) and (324.94,113.74) .. (327.67,101) .. controls (330.88,86) and (334.56,69.77) .. (335.67,54.33) ;
\draw  [color={rgb, 255:red, 155; green, 155; blue, 155 }  ,draw opacity=1 ][line width=0.75] [line join = round][line cap = round] (385,73.67) .. controls (385.73,72.93) and (392.24,49.26) .. (395,45) .. controls (396.11,43.28) and (398.89,47.62) .. (399,49.67) .. controls (399.97,68.42) and (397.67,85.02) .. (397.67,103) ;
\draw  [color={rgb, 255:red, 155; green, 155; blue, 155 }  ,draw opacity=1 ][line width=0.75] [line join = round][line cap = round] (438.97,181.83) .. controls (439.29,190.58) and (449.71,229.7) .. (462.13,218.38) .. controls (466.32,214.56) and (468.4,206.73) .. (470.12,201.64) ;
\draw    (159,119) -- (205,95) ;
\draw  [dash pattern={on 4.5pt off 4.5pt}]  (160.72,69.31) -- (409.67,213) ;
\draw  [fill={rgb, 255:red, 245; green, 166; blue, 35 }  ,fill opacity=0.6 ] (205,95) -- (160.9,118.46) -- (160.74,69.59) -- cycle ;
\draw    (178.67,92) -- (223.04,39.86) ;
\draw [shift={(224.33,38.33)}, rotate = 130.4] [color={rgb, 255:red, 0; green, 0; blue, 0 }  ][line width=0.75]    (10.93,-3.29) .. controls (6.95,-1.4) and (3.31,-0.3) .. (0,0) .. controls (3.31,0.3) and (6.95,1.4) .. (10.93,3.29)   ;
\draw    (357.33,190.67) -- (328.57,221.08) ;
\draw [shift={(327.2,222.53)}, rotate = 313.4] [color={rgb, 255:red, 0; green, 0; blue, 0 }  ][line width=0.75]    (10.93,-3.29) .. controls (6.95,-1.4) and (3.31,-0.3) .. (0,0) .. controls (3.31,0.3) and (6.95,1.4) .. (10.93,3.29)   ;
\draw    (406,48) -- (432.65,30.75) ;
\draw [shift={(434.33,29.67)}, rotate = 147.09] [color={rgb, 255:red, 0; green, 0; blue, 0 }  ][line width=0.75]    (10.93,-3.29) .. controls (6.95,-1.4) and (3.31,-0.3) .. (0,0) .. controls (3.31,0.3) and (6.95,1.4) .. (10.93,3.29)   ;
\draw [color={rgb, 255:red, 155; green, 155; blue, 155 }  ,draw opacity=1 ]   (143,117.67) .. controls (155,101.67) and (157.67,76.33) .. (163.67,61) .. controls (169.67,45.67) and (171,23.12) .. (181,21.67) .. controls (191,20.22) and (197,122.33) .. (196.33,135.67) ;

\draw (509.33,137.07) node [anchor=north west][inner sep=0.75pt]    {$\delta _{i}$};
\draw (464.33,80.4) node [anchor=north west][inner sep=0.75pt]  
{$\frac{1}{3} \delta _{i}$ };
\draw (137.67,183.07) node [anchor=north west][inner sep=0.75pt]  
 {$\frac{1}{3} \delta _{i}$ ~};
\draw (465.33,130.4) node [anchor=north west][inner sep=0.75pt] 
 {$\frac{1}{3} \delta _{i}$};
\draw (416,125.4) node [anchor=north west][inner sep=0.75pt]  
 {$\widetilde{u}_{k}$};
\draw (107.7,65.4) node [anchor=north west][inner sep=0.75pt]  
{$v_{i}{}_{-1}$};
\draw (447.7,228.4) node [anchor=north west][inner sep=0.75pt] 
{$v_{i}$};
\draw (226.67,25.73) node [anchor=north west][inner sep=0.75pt]  
 {$T_{i}$};
\draw (297.33,232.4) node [anchor=north west][inner sep=0.75pt] 
 {$\mathrm{slope} \ = \  -2^{k}$};
\draw (440.67,16.07) node [anchor=north west][inner sep=0.75pt]  
  {$\overline{u}$};

\end{tikzpicture}

    \caption{\small The estimate for $i \in {\cal J}_2$. From the fact that the value $v_{i-1}$ survives up to time $2^{-k}$, we deduce that the area of the yellow triangle can be controlled by the $L^1$ norm of all the $v_q^p$ in the interval $(x_{i-1}, x_i)$.}
    \label{fig:mainest}
\end{figure}

The proof splits in two cases depending on which one is true. 
To fix ideas, we assume that the first inequality holds.  The second case is entirely similar. It can be handled  by the same argument, in connection with the reversed initial datum:
$
\ol v(x) \doteq -\ol u(-x)
$.

Assuming that the first inequality in (\ref{in12}) holds,
define the triangle:
$$
\T_i ~\doteq ~\Big\{ (x, v) \in \mathbb R^2 \; \mid \; x \in (x_{i-1}, x_i), \quad v_{i-1} -\frac{\delta_i}{3} + (x-x_{i-1})\cdot 2^{k} \leq v \leq v_{i-1}-(x-x_{i-1})2^k \Big\},
$$
as shown in Fig.~\ref{fig:mainest}.
We let $U \subset \mathbb R^2$ be the hypograph of $\ol u$:
$$
U~ \doteq ~\Big\{ (x, v) \in \mathbb R^2 \,;~~ v \leq \ol u(x) \Big\}~ \subset ~\mathbb R^2,
$$
and let $\wt U_k$  be the hypograph of $\wt u_k$:
$$
\wt U_k~ \doteq ~\Big\{ (x, v) \in \mathbb R^2 \,;~~  v \leq \wt u_k(x) \Big\} 
~\subset~ \mathbb R^2.
$$

The fact that  the couple $( x_{i-1},\,v_{i-1})$ survives up to time $t=2^{-k}$ 
already implies that 
\begin{equation}\label{eq:survarea}
\meas(\T_i\setminus \wt U_k) ~\leq ~\meas \Big\{ (x,v)\in (U \setminus \wt U) \,;~~x\in [x_{i-1}, x_i] \Big\}~\doteq~A_i\,.
\end{equation}

Actually, we claim that $\wt U_k \cap \T_i = \emptyset$. In fact, $\wt u_k$ is one-sided $2^k$-Lipschitz and satisfies $\wt u_k(x_{i-1}) \leq v_{i-1}-\frac{\delta_i}{3}$. This
implies
$$
\wt u_k(x)~ \leq~ v_{i}-\frac{\delta_i}{3}+ (x-x_{i-1})\cdot 2^k\qquad \forall \; x \in 
\,]x_{i-1}, x_i[\,.
$$
By definition of the triangle $\T_i$, this  means that the hypograph of $\wt u_k$ lies entirely below the lower side of $\T_i$. Hence $\wt U_k \cap \T_i = \emptyset$.  From \eqref{eq:survarea} it thus follows
\begin{equation}\label{eq:survarea1}
\meas(\T_i)~ \leq ~A_i\,.
\end{equation}
On the other hand, the area of the triangle $\T_i$ is 
\begin{equation}\label{eq:triarea}
\meas(\T_i)~ =~ \frac{\delta_i^2}{64} \,2^{-k}.
\end{equation}
Combining \eqref{eq:survarea1} with  \eqref{eq:triarea} we obtain
\begin{equation}\label{eq:triest1}
\delta_i ~\leq ~8 \cdot 2^{k/2}\cdot A_i^{1/2}.
\end{equation}
The area $A_i$ on the right hand side is bounded above by the sum of the areas  of the blocks $v_q^p$ whose support intersects the interval $\,]x_i, x_{i+1}[\,$.
More precisely, for each $i \in \mathcal{J}_2$, define the set of indices 
$$
{\cal Z}(i) \doteq \Big\{ (p, q) \in \mathbb N^2 \; \mid \; q \geq k, \quad {\rm supp}\, v_q^p \cap \,]x_{i-1}, x_i[\,  \neq \,\emptyset\Big\}.
$$
With the same argument used in step {\bf 5.} we  obtain that,
for every couple
$(p, q)$  with $q \geq k$, there can be at most two indices $ i \in {\cal J}_2$ 
such that  $ (p, q) \in {\cal Z}(i)$.


By \eqref{eq:triest1} we now obtain
$$
\delta_i ~ \leq ~ 8\cdot 2^{k/2}\cdot \Bigg(  \,\sum_{(p,q) \in {\cal Z}(i)} h_q^p \cdot \ell_q^p\, \Bigg)^{1/2}~ \leq~  8\cdot 2^{k/2}\cdot  \sum_{(p,q) \in {\cal Z}(i)} \Big( h_q^p \cdot \ell_q^p\, \Big)^{1/2} .
$$
Summing  over $i$, and using the fact that each $(p,q)$ appears in the sum at most twice,  we obtain 
\begin{equation}\label{eq:twoiineq1}
    \sum_{i \in {\cal J}_2} \delta_i ~ \leq~ 8 \cdot 2^{k/2}\cdot  \sum_{i \in {\cal J}_2} \, \sum_{(p,q) \in {\cal Z}(i)} \Big( h_q^p \cdot \ell_q^p\, \Big)^{1/2} ~ \leq ~ 16 \cdot 2^{k/2} \cdot \sum_{q = k}^{\infty}\, \sum_{p\in \mathbb N}\Big( h_q^p \cdot \ell_q^p\, \Big)^{1/2}.
\end{equation}
Observing that 
$
\Big(\ell_q^p\cdot  h_q^p\Big)^{1/2} = 2^{q/2}\ell_q^p
$
and using \eqref{eq:rectsum}, we finally obtain
\bel{es63}\bega{rl}\ds
\sum_{i \in {\cal J}_2} \delta_i & \leq\ds~   16 \cdot 2^{k/2} \cdot \sum_{q=  k}^{\infty} \, \sum_{p \in \mathbb N}  \, 2^{q/2}\ell_q^p ~\leq~16 \cdot 2^{k/2} \cdot  C\cdot \sum_{q =k}^{\infty}\,  2^{q/2}  2^{-\alpha q}\\[4mm]
 & \leq ~\ds c_3 \cdot \frac{1}{2\alpha-1}\cdot   \Vert \ol u\Vert_{{\cal P}_{\alpha}}2^{k/2} \cdot 2^{(1/2-\alpha) k} ~=~ c_3 \cdot \frac{1}{2\alpha-1}\cdot \Vert \ol u\Vert_{{\cal P}_{\alpha}} \cdot 2^{(1-\alpha) k},
\enda
\eeq
where $c_3$ is an absolute constant. 
Combining the three estimates (\ref{es61}), (\ref{es62}) and (\ref{es63}), 
the proof is completed.
\endproof
\v

\begin{remark}\label{r:62}  {\rm
 If Burgers' equation is replaced by a general scalar conservation law with a $\C^2$,
 uniformly convex flux $f$, so that $f^{\prime \prime } \geq c  > 0$, from the Hopf-Lax formula we obtain that $(x_0, u_0) \in \Q(t)$  if and only if 
 \begin{equation}\label{eq:rightsurvgen}
 \int_{x_0}^y \left[\ol u(z)-(f^*)^{\prime}\Big(f^{\prime}(u_0)-\frac{z-x_0}{t} \Big)\right]^+ \,d z ~
 \geq~  \int_{x_0}^y \left[\ol u(z)-(f^*)^{\prime}\Big(f^{\prime}(u_0)-\frac{z-x_0}{t} \Big)\right]^- d z   
        \end{equation}
  for all    $ y \geq x_0\,$, and
    \begin{equation}\label{eq:leftsurvgen}
         \int_{y}^{x_0} \left[\ol u(z)-(f^*)^{\prime}\Big(f^{\prime}(u_0)-\frac{z-x_0}{t} \Big)\right]^+ d z ~\leq ~ \int_{y}^{x_0} \left[\ol u(z)-(f^*)^{\prime}\Big(f^{\prime}(u_0)-\frac{z-x_0}{t} \Big)\right]^- d z \eeq
         for all  $y \leq x_0\,$.
Here 
$$f^*(u)~\doteq~\sup_{v\in \R} ~\bigl\{uv-f(v)\bigr\}$$ denotes 
the Legendre transform of $f$.  
By a well known property of the Legendre transform (see e.g. \cite{Canna_book, E}) we have
$$
(f^*)^{\prime}\bigl(f^{\prime}(u_0)\bigr) ~=~ u_0\,.
$$
By uniform convexity it thus follows
 $$
 u_0 -\frac{1}{ct}(z-x_0)~\leq ~(f^*)^{\prime}\Big(f^{\prime}(u_0)-\frac{z-x_0}{t}\Big)  \qquad \forall \; z > x_0\,,
 $$
$$
(f^*)^{\prime}\Big(f^{\prime}(u_0)-\frac{z-x_0}{t}\Big) ~ \leq ~u_0 -\frac{1}{ct}(z-x_0) \qquad \forall \; z < x_0.
$$
Using the inequalities above, the above proof remains valid up to minor modifications. Therefore, Theorem \ref{thm:dec} remains valid for any uniformly convex flux $f$.

}
\end{remark}

{\bf Acknowledgements.}  The research of A.~Bressan
 was partially supported by NSF with
grant  
  DMS-2306926, ``Regularity and approximation of solutions to conservation laws".
F.~Ancona, E.~Marconi and L.~Talamini are members of GNAMPA of the ``Istituto Nazionale di Alta Matematica
F.~Severi". F.~Ancona and E.~Marconi are partially supported by the PRIN 2020 project “Nonlinear evolution PDEs, fluid dynamics and transport
equations: theoretical foundations and applications”.
E.~Marconi  is also partially supported by H2020-MSCA-IF “A~Lagrangian approach: from conservation laws to line-energy Ginzburg-Landau models”.

\end{document}